\documentclass{amsart}
\usepackage{amsfonts,amscd}

\newtheorem{claim}{}[section]
\newtheorem{theorem}[claim]{Theorem}
\newtheorem{lemma}[claim]{Lemma}
\newtheorem{proposition}[claim]{Proposition}
\newtheorem{corollary}[claim]{Corollary}

\newtheorem{definition}[claim]{Definition}

\def\proclaim #1. #2\par{\medbreak
\noindent{\bf#1.\enspace}{\sl#2}\par\medbreak}
\makeatother

\DeclareMathOperator{\Ndb}{\mathbb{N}}

\DeclareMathOperator{\Tdb}{\mathbb{T}} 
\DeclareMathOperator{\Rdb}{\mathbb{R}}
\DeclareMathOperator{\fM}{\mathfrak{M}}   
\DeclareMathOperator{\cas}{\text{C*-algebras}}

\DeclareMathOperator{\sub}{\subseteq}

\begin{document}

\title[Ordered C*-modules]
{Ordered C*-modules}

\thanks{Some of the contents of this paper was summarized at the 
G.P.O.T.S. conference at the University of
Urbana-Champaign, May 2003.}
\thanks{This research was supported in part by the
SFB 487  Geometrische Strukturen in der Mathematik,
at the Westf\"alische Wilhelms-Universit\"at, supported
by the Deutsche Forschungsgemeinschaft.
Blecher was also supported in part by a grant from the
National Science Foundation.} 

 \author{David P. Blecher}
\address{Department of Mathematics, University of Houston, Houston,
TX
77204-3008, USA}
\email[David P. Blecher]{dblecher@math.uh.edu}

\author{Wend Werner}
\address{SFB 478 - Geometrische Strukturen in der Mathematik,
Westfalische Wilhelms-Universitat, Hittorfstrasse 27, 48149 Munster, Germany}
\email[Wend Werner]{wwerner@math.uni-muenster.de}

\maketitle

\let\text=\mbox

\begin{abstract}   In this first part of a study of
ordered operator spaces, we develop the basic theory of  
`ordered C*-bimodules'.  A crucial role is played by
`open central tripotents', a JB*-triple variant of 
Akemann's notion of open projection.
\end{abstract}

\section{Introduction}

In \cite{HBWI,Weaver} Weaver introduced Hilbert 
C*-bimodules with involution.  In the present paper we consider
a further structure on C*-bimodules: order.  That is, we study
C*-bimodules which have an involution and a positive cone.
The reason why we became interested in this topic is
that such C*-bimodules arise naturally as the `noncommutative
Shilov boundaries' of ordered operator spaces (the latter term is 
defined below).  
Two consequences of this: firstly, ordered C*-modules are excellent 
and appealing examples
of ordered operator spaces; and secondly, general  
ordered operator spaces may be studied using  ordered C*-modules.
We develop in the present paper 
the basic theory of ordered C*-bimodules, or what is
equivalent, ordered {\em ternary rings of operators} (defined
below).  This theory is in some sense quite uncomplicated
(for example the only `ordered W*-bimodules' are the sum
of a W*-algebra and an `unorderable W*-bimodule'), but in other 
ways does possess intricate features.
As an interesting byproduct,
it turns out that a number of properties of C*-bimodules that, in
general, are
only definable using an embedding into a C*-algebra actually turn out to be
invariants of the underlying order structure.   We introduce some of
these invariants in Sections 2 and 3, where we principally discuss involutive
 structure. 
In 
Section 4, we characterize the possible orderings on a C*-bimodule $X$ which
correspond to embeddings of $X$ as selfadjoint
 ternary rings of operators.   Our characterization is   
in terms of selfadjoint {\em tripotents};
namely elements $u$ such that $u = u^* = u^3$.  In the uniform version
of our theory, we will need the `tripotent' variant of Akemann's notion 
(see e.g.\ \cite{Ake}, \cite[3.11.10]{Ped}) of an open projection.
Hence we will also study in Section 4 
some of the basic properties associated with this key notion.
In Section 5 we characterize all the {\em maximal} ordered operator space
orderings on $X$, these turn out to be 
automatically orderings of the type discussed above.
In Section \ref{ilb} 
we consider the important `commutative' example, 
{\em involutive line bundles}.  Amongst other things
these will show that our results in previous sections are sharp.
In the final section we study some relations with a certain topic
in physics.  In a sequel paper we will apply some of our present results
to analyze general ordered operator spaces, using their
noncommutative
Shilov boundaries mentioned above.  

We now turn to precise definitions and notation.  
We write $X_+$ for the cone of `positive elements',
i.e.\ those with $x \geq 0$, in an ordered vector space.
A subspace $X$ of a vector space with involution is {\em selfadjoint} 
if $x^* \in X$ whenever $x \in X$,  or in shorthand,
$X^* = X$.  We will say that a linear map $T : X \to Y$ between
vector spaces with involution is {\em selfadjoint}
or {\em $*$-linear} if $T(x^*) = T(x)^*$ for all $x \in X$.
If $X$ and $Y$ are also ordered vector spaces then
we say that $T$ is {\em positive} if
$T$  is selfadjoint, and if   $T(X_+) \subset Y_+$. 
By a {\em concrete ordered operator space}, we mean a
selfadjoint subspace $X$ of a C*-algebra $A$, 
together with the positive cone $X \cap A_+$ inherited from $A$. 
If also $X$ contains
$1_A$ we call $X$ a {\em unital operator system}.    
We will not use much from operator space theory, what we use can 
probably be found in any of the current books, or in most of the papers, on
operator spaces.  We use the term {\em ordered operator space}
for an operator space with an involution,
and with specified positive cones in $M_n(X)_{sa}$, so that there exists a
completely isometric complete order embedding from $X$ into a C*-algebra. 
We will not need this here, but such objects were characterized abstractly
in \cite{OSWU}.  We say that one
ordering on $X$ is {\em majorized} by another ordering if the positive cones
for the first ordering are contained in the positive cones for the second
ordering.

We have used $*$ above for the `adjoint' or `involution';
because this symbol appears so frequently in this paper we will 
instead write $X'$ for the dual Banach space
(resp.\ dual operator space) of a Banach space
(resp.\ operator space) $X$.
 
By a {\em ternary ring of
operators} (or {\em TRO} for short) we will mean a closed subspace
$X$ of a C*-algebra $B$ such that $X X^* X \subset X$.
Sometimes we call this a TRO {\em in } $B$.
The important structure on a TRO is the `ternary product'
$x y^* z$, which we sometimes write as 
$[x,y,z]$.   Our interest here is in TRO's
within the realm of operator spaces,
i.e.\  we will always regard a TRO $X$ as equipped with its canonical
operator space structure.  Thus there are canonical
norms on $M_n(X)$ for each $n \in \Ndb$
(which are compatible with the
natural TRO structure on $M_n(X) \subset M_n(B)$).
Under this assumption, an intrinsic abstract
characterization of TRO's may be found in
\cite{NeRu} (see also \cite{ACOT}).  We will not use this
however.    Instead we simply make a definition: a {\em 
ternary system} is an operator space $X$ with a map
$[\cdot,\cdot,\cdot] : X \times X \times X \to X$, such that  
there exists a TRO $Z$ and a
completely isometric surjective linear isomorphism $T : X \to Z$ satisfying 
$[Tx,Ty,Tz] = T([x,y,z])$ for all $x, y, z \in X$.
A linear map $T$ satisfying the last equation is called a 
{\em ternary morphism}.  Even for an
abstract ternary system we sometimes simply write
 $[x,y,z]$ as $x y^* z$ for $x, y, z \in X$.

We write $X Y$ for the norm closure of the span of the set
of products $x y$ for $x \in X, y \in Y$, assuming that such
products make sense.  A similar notation applies to the 
product of three or more sets.  In contrast, in the
expression $X + Y$
we are not automatically taking the norm closure.
It is well known  that $X X^* X = X$
for a TRO $X$.  Also,
it is clear that $X X^*$ and $X^* X$ are C*-algebras,
which we will call the {\em left} and 
{\em right C*-algebras} of $X$ respectively,
and $X$ is a $(X X^*)$-$(X^* X)$-bimodule.
 
The following result shows that
ternary morphisms behave very similarly to *-homomorphisms
between C*-algebras:

\begin{proposition} \label{prtm} {\rm (See e.g.\ \cite{Ham3})}
Let $T : X \to Y$ be a ternary morphism between ternary systems.  Then:
\begin{itemize} \item [(1)]  
$T$ is completely contractive and has closed range.
\item [(2)]  $T$  is
completely isometric if it is 1-1.  
\item [(3)]  A linear isomorphism between ternary systems is
completely isometric if and only if it is a ternary morphism.  
\item [(4)]  The quotient of a TRO $X$ by a `ternary ideal' (that is, a
uniformly closed $(X X^*)$-$(X^* X)$-subbimodule)
is a ternary system.
\item [(5)]   If $X$ is a TRO 
then {\rm Ker}$(T)$ is a ternary ideal,
and the induced map
$X/{\rm Ker}(T) \to Y$ is a 1-1 ternary morphism.
\item [(6)]  If $X, Y$ are TRO's,
 then $T$
canonically induces a $*$-homomorphism (resp.  $*$-isomorphism)
$\pi : X^* X
\rightarrow Y^* Y$ between the associated
right C*-algebras, via the prescription $\pi(x^* y) = T(x)^* T(y)$.
If $T$ is 1-1 (resp.
1-1 and onto) then $\pi$ is also 1-1 (resp.
a $*$-isomorphism).  
Similarly for the left C*-algebras. 
\end{itemize} \end{proposition}

From (3) of the last Proposition one easily sees that 
an operator space $X$ may have at most one
`ternary product' $[\cdot,\cdot,\cdot]$ with respect
to which it is a ternary system.   Thus we may simply define a
ternary system to be an operator space $X$ which is linearly
completely isometric to a TRO.

In this paper we will be more interested in {\em $*$-TRO's},
by which we mean
a closed selfadjoint TRO $Z$ in a C*-algebra $B$.
In this case,
the `left' and `right' C*-algebras mentioned above,
namely $Z Z^*$ and $Z^* Z$,
coincide, and equal $Z^2$.
By a {\em $*$-WTRO} we mean a selfadjoint weak*
closed TRO in a W*-algebra.  Note that a $*$-TRO $Z$ comes
 with a given positive cone $Z_+$, inherited 
from the containing  C*-algebra.      
Analogues of parts of the last proposition are 
valid for selfadjoint ternary morphisms,
which we also call  {\em ternary $*$-morphisms},
between $*$-TRO's.  In particular, 
the kernel of a ternary $*$-morphism on a $*$-TRO
is clearly a  {\em ternary $*$-ideal}, by which we mean
a  selfadjoint ternary ideal.
It is not hard, by following the usual
proof of Proposition \ref{prtm} (4) (or by using 
Lemma \ref{*-TRO_char} (1) below),
to show that the quotient
of a $*$-TRO by a ternary $*$-ideal is
again an involutive ternary system in a natural way which is
compatible with the quotient operator space structure on $Y/X$.
From this one easily checks the analogue of 
 Proposition \ref{prtm} (5); that if one factors a ternary $*$-morphism
on a $*$-TRO by
its kernel, then one obtains a 1-1 ternary $*$-morphism on the
quotient ternary system. 
A ternary $*$-ideal $J$ in an involutive ternary system $Z$
will be called a {\em C*-ideal} if
$J$ is ternary $*$-isomorphic to a C*-algebra.  

Ternary systems, and hence also $*$-TRO's, are a particularly nice 
subclass of the {\em JB*-triples}.   We will not define the 
latter objects, but we  must stress that many of the techniques and 
ideas in the present paper originate in that field of study.
As we said above, we establish a link between 
orderings on a $*$-TRO $Z$, and certain `central' selfadjoint
tripotents in $Z$ or $Z''$.
With this in hand,
a certain portion of our results may be viewed as variants of certain
JB*- and JBW*-triple results.  In spite of this we have avoided 
deriving many such results of
ours from appeals to the JBW*-triple literature, for several reasons.
Firstly, the proofs of which we are speaking are quite short and simple
in our setting, and so it seemed much more natural to include the direct
arguments.  Secondly, we could not find such results
in the precise form we needed, and if we had to include 
the details of the modifications of
the JBW*-triple results, the proofs would become unnecessarily long.  
Then of course many of the 
issues in the JB*-triple theory are not relevant in our setting (being automatic
or simply do not arise).  However we have tried to consistently
indicate, to the best of our knowledge,
where a comparison with the JBW*-triple literature should be 
made.   
The reader should  certainly
browse concurrently with our paper the
JB*-triple papers; for example \cite{BaTim,Bat,FRc,AGOC,Horn2,IKR,NeRu2,SOJ};
and also the mammoth edifices of work by W. Kaup, and 
by C. M. Edwards and the late G. T. R\"uttimann, of which we have cited some 
representatives which have some important points of contact 
with our paper.  Note that if $Z$ is a $*$-TRO then $Z_{sa}$ is a real 
JB*-triple (see e.g.\ \cite{ER4,IKR,SOJ}), and one may then
appeal to the methods and results of the 
real JB*-triple theory.  As far as we 
are aware however, the main results of our paper are quite new.
In particular, we have not seen `open tripotents' in our sense in
the literature.
(In \cite{ER3,ER2} this term is used in a sense which is formally related to ours,
but is different.  
Also there is much work on `compact tripotents' (see e.g.\ \cite{ER2}), which in 
any case is not what is needed here.  These works
were partially inspired by Akemann and Pedersen's paper \cite{AkeP}).
 
In a similar spirit, we should say that we use `TRO techniques'
throughout; and  there are several recent papers of 
Ruan alone or with coauthors
concerning TRO's (see e.g.\ \cite{OIAN,Ruan}).  Again
any overlaps between this work and ours, only concerns very simple facts.
  In passing we remark that Proposition \ref{prtm} (3) is independently 
due to Ruan.
 
We end this section with a note on the title of this paper,
{\em ordered C*-modules}.  An ordered C*-module
 is a C*-bimodule $Y$ over a C*-algebra $A$
with a given involution and positive cone, such that
$Y$ (with its canonical ternary product $x \langle y , z \rangle$)
is `ternary order isomorphic' to  a $*$-TRO.
Thus ordered C*-modules are essentially the same thing
as $*$-TRO's in a C*-algebra $B$, 
with their canonical inherited ordering from $B$.  
In the last paragraph of Section \ref{order},
we will give a better  characterization 
of ordered C*-modules 
amongst the `involutive ternary systems' (defined below).
Since they are essentially the same as $*$-TRO's, it
suffices to focus on the order properties of $*$-TRO's:
nearly all of our results on $*$-TRO's will 
transfer in an obvious way to ordered C*-modules.
Thus the reader will not see the term ordered C*-module
much in this paper.

\smallskip

{\em Acknowledgements.}   Both authors are 
very grateful for support from the SFB 487  Geometrische Strukturen in der Mathematik,
at the Westf\"alische Wilhelms-Universit\"at, in turn supported
by the Deutsche Forschungsgemeinschaft.
The first author is also indebted 
to our colleagues there for their generous hospitality,
in 2002 when this work began.  
Another impetus for this project was
 a Research Experiences for Undergraduates problem that
the first author gave to then undergraduate Kay Kirkpatrick \cite{TSBA}.
Finally we thank T. Oikhberg, and B. Russo 
for several comments after a talk the first author gave on this work;
and M. Neal for much helpful information.  
  
\section{Involutive $C^*$-bimodules} \label{ICB}

In this section we fix a C*-algebra $B$, which the
reader may wish to take equal to $B(H)$ for a Hilbert space $H$.  Then
$M_2(B)$ is also canonically a C*-algebra.  Indeed
$M_2(B) \cong B(H^{(2)})$ $*$-isomorphically, if
$B = B(H)$.  Let  $A$ and $Z$ be selfadjoint subspaces of $B$,
and define 
\begin{equation} \label{bzzb}
{\mathcal L} \; = \; \left[
\begin{array}{cc}
A & Z \\
Z & A
\end{array}
\right]  \; = \; \left\{  \left[
\begin{array}{cc} a_1 & z_1 \\ z_2 & a_2 \end{array}
\right] \in M_2(B) : a_1, a_2 \in A, z_1, z_2 \in Z \right\} .
\end{equation} 
Then ${\mathcal L}$ is a selfadjoint 
subspace of $M_2(B)$, and it  
is uniformly closed if $A$ and $Z$ are closed 
 in $B$. 
Similarly, if $B$ is a W*-algebra, and if 
$A$ and $Z$ are weak* closed in $B$, then
${\mathcal L}$
is  weak* closed in the W*-algebra $M_2(B)$.

Similar remarks apply to the selfadjoint closed subspace 
\begin{equation} \label{bzzb2}
\tilde{{\mathcal L}} = \left\{ \left[
\begin{array}{cc}
a & z \\
z & a \end{array}
\right]  \in {\mathcal L} : a \in A, z \in Z \right\} \end{equation}
of
${\mathcal L}$.   
Note also that $\tilde{{\mathcal L}}$
is canonically completely isometric 
to a subspace  of $B \oplus^\infty B$.
To see this, notice that the maps
$$\left[
\begin{array}{cc}
x & y \\
y & x 
\end{array}
\right] \mapsto \left( x + y , x - y \right),$$
and $$ \left( x , y  \right)  \mapsto  
\left[ \begin{array}{cc}
\frac{x+y}{2} & \frac{x-y}{2} \\ 
\frac{x-y}{2} &  \frac{x+y}{2}  
\end{array}
\right],
$$
are $*$-isomorphisms between $B \oplus^\infty B$
and the subspace
$$
\{ \left[
\begin{array}{cc}
x & y \\
y & x
\end{array}
\right] : x , y \in B \}
$$
of $M_2(B)$.  
In particular it follows that 
\begin{equation} \label{xyyxnorm}  
\left| \left| \left[
\begin{array}{cc}
x & y \\
y & x
\end{array}
\right] \right| \right| = 
\max \{ \Vert x + y \Vert ,
\Vert x - y \Vert \} .
 \end{equation}

Henceforth we will consider the case when ${\mathcal L}$ 
and $\tilde{{\mathcal L}}$ above are C*-subalgebras of $M_2(B)$.
This occurs precisely when $A$ is a C*-subalgebra
of  $B$, and $Z$ is a closed selfadjoint subspace
of $B$ such that $AZ \subset Z$ and $Z^2 \subset A$.
In this case we say that $Z$ is an {\em involutive ternary $A$-submodule
of} $B$.   Then $\tilde{{\mathcal L}}$ above
is canonically  $*$-isomorphic to a $*$-subalgebra
 of $B \oplus^\infty B$.  
 We note that any involutive ternary $A$-submodule is 
a $*$-TRO.  Conversely, any $*$-TRO $Z$ in $B$ is an
involutive ternary $Z^2$-submodule
of $B$.    
If $B$ is a  W*-algebra, and if $Z$ and $A$ are weak* closed
in $B$, then $Z$ is a $*$-WTRO, and
also  ${\mathcal L}$
and $\tilde{{\mathcal L}}$  are W*-subalgebras of $M_2(B)$.

There are abstract characterization of involutive ternary $A$-bimodules.
These bimodules are in fact very closely related to the 
topic of $*$-automorphisms $\theta$ of a C*-algebra $M$
of {\em period 2}, that is, $\theta^2 = Id$.
  (Although 
we shall not use this, we observe in passing the obvious fact
that such automorphisms are in a bijective correspondence 
with involutions on $M$, commuting with 
the usual involution, and with respect to which $M$ is still
a C*-algebra.)    
For such an automorphism  
$\theta$, let $N$ be the fixed point algebra $\{ a \in M :
\theta(a) = a \}$, and let $W = \{ x \in M : \theta(x) = - x \}$.
Clearly $W$ is an involutive ternary $N$-subbimodule of $M$.
Conversely, if $Z$ is an involutive ternary $A$-subbimodule 
in $B$, then the C*-algebra  $\tilde{{\mathcal L}}$ in Equation
(\ref{bzzb2})
has an obvious
period 2 automorphism whose fixed point algebra
is isomorphic to $A$; and $Z$ is appropriately 
isomorphic to the set of matrices $x \in \tilde{{\mathcal L}}$ 
with $\theta(x) = - x$.
Indeed if $Z \cap A = (0)$, then one can show that
$M = Z + A$ is a C*-algebra with an obvious
period 2 automorphism $z + a \mapsto z - a$ (see 
Corollary \ref{iscl} (1) below).   

The following abstract characterization of involutive ternary $A$-bimodules,
while conceptually significant,
will not be technically important in the present paper,
so we merely take a couple of paragraphs to sketch the details:
 
\begin{definition} \label{incb}
An {\em involutive C*-bimodule} over
a C*-algebra $A$ is a bimodule
$X$ over $A$,  such that $X$ has an involution
satisfying $(ax)^* = x^* a^*$, and such that $X$
is both  a  right and a  left
C*-module over $A$, with the left module inner product
$[\cdot|\cdot]$ being related to the right
module inner product by the formula $[x|y] = \langle x^*|y^* \rangle$,
and such that two inner products are also
compatible in the sense that
$x \; \langle y  \vert   z \rangle  =
[ x   \vert   y ] \; z$ for all $x,y,z \in X$.
 
We say that such a bimodule $X$ is a {\em
commutative involutive C*-bimodule} if in addition
$\langle y  \vert   z \rangle = [z,y]$ for all $y,z \in X$.
\end{definition}
 
Since an involutive C*-bimodule is a
C*-bimodule, it has a canonical operator space
structure (see e.g.\ \cite{Blecher}).
  
It is clear that any  involutive ternary $A$-submodule of a
C*-algebra $B$, is an involutive C*-bimodule in the above sense.
To see the converse, we note
that the proof of Theorem 12 in \cite{HBWI}
may be easily adapted to show that any
involutive C*-bimodule over a C*-algebra $A$ is isomorphic,
 via a complete isometry which is also a
ternary $*$-isomorphism and an $A$-$A$-bimodule map,
to an involutive ternary $A$-submodule of a C*-algebra $C$.
See also our later Corollary \ref{wcor}.
  
We will not use this,
but involutive C*-bimodules may be equivalently defined to
be a certain subclass of
the Hilbert $A$-bimodules with
involution recently developed by Weaver.
  Namely, it is the subclass of bimodules satisfying the
natural extra condition $( x , y ) z = x ( y , z )$
in the language of Weaver's definition on the second page
of \cite{HBWI}, and
which are also Hilbert C*-modules in the usual sense with regard
to both of the natural sesquilinear $A$-valued inner
products.

\medskip
              
There are three important spaces canonically associated to 
a $*$-TRO, or involutive ternary $A$-submodule, $Z$ in a
C*-algebra $B$.
The first is the $*$-subalgebra $A + Z$ of $B$, which we will see 
momentarily  is a C*-algebra.   The second space
 we call the
{\em center} of $Z$ (following
\cite{HBWI,Sau,Skei}): this is defined to be
$${\mathcal Z}(Z) \; = \; \{ z \in Z : a z = z a  \;
\text{for all} \; a \in Z^2 \}. $$    
It is easy to see that ${\mathcal Z}(Z)$ is a 
$*$-TRO in $B$ too.  One can show that  
$z \in {\mathcal Z}(Z)$ if and only if $a z = z a$ for all 
$a \in A$, but since we will not use this fact
we will not prove it here.
We will study ${\mathcal Z}(Z)$ in more detail later in this 
Section.  The third space $J(Z)$ is defined to be $Z \cap Z^2$. 
Clearly we have
$$J(Z) \subset Z \subset Z + A .$$
 
These three auxiliary spaces
will play a significant role for us.   Note 
that ${\mathcal Z}(Z)$ is an invariant of the involutive ternary
structure on $Z$.  That is,
a ternary $*$-isomorphism $\psi : E \to F$ between two
$*$-TRO's, restricts to a
ternary $*$-isomorphism from ${\mathcal Z}(E)$ onto ${\mathcal Z}(F)$.
This may be easily seen using Proposition 
\ref{prtm} (6).   On the other hand, 
$Z^2 + Z$ and $J(Z)$ 
are not invariants of the involutive  ternary 
structure on $Z$.
However we shall see later in Corollary \ref{contof},
the interesting fact that
$Z^2 + Z$ and $J(Z)$ {\em are}  invariants of the 
order structure of $Z$.   For this reason we defer much
discussion of $Z^2 + Z$ and $J(Z)$ to the later Section
\ref{order} which is devoted to order structure.

\begin{corollary} \label{iscl}  Let  $Z$ be an
involutive ternary $A$-submodule of $B$.
Then
\begin{itemize}
\item [(1)]  $A + Z$ is a closed subspace
of $B$, and indeed is a
C*-subalgebra of $B$.  
\item [(2)]  If further $B$ is a W*-algebra, and if $A$ and $Z$ are
weak* closed in $B$, then $A + Z$ is a weak* closed subspace
of $B$, indeed is a W*-subalgebra of $B$.   
\item [(3)]  $J(Z) = Z \cap A$.
\item [(4)]   $J(Z)$ is an ideal in both
$A$ and $A + Z$.  
\item [(5)]  $J(Z)$ is a complete $M$-ideal in $Z$. 
It is also an $A$-subbimodule of $Z$,
and a ternary $*$-ideal in $Z$.
\end{itemize}
\end{corollary} 

\begin{proof}  
(1) \ As remarked above, 
in this case the set $\tilde{{\mathcal L}}$
in Equation (\ref{bzzb2})
is a C*-algebra.   The map $\pi$ from 
$\tilde{{\mathcal L}}$ to $B$ taking the matrix
written in Equation (\ref{bzzb2}) to $a + z$, is a $*$-homomorphism.  
Therefore its range is a closed C*-subalgebra of $B$.

(2)  \ Assuming the hypotheses here,
then as remarked above, the set $\tilde{{\mathcal L}}$
is a W*-algebra.  It is also easy to see that 
$\pi$ in the proof of (1)
is weak* continuous.  By basic von Neumann algebra
theory we deduce that Ran$(\pi)$ is weak* closed.

(3) \ Suppose that $z \in Z \cap A$.  If $(e_t)$ is an increasing
approximate identity for $Z^2$, then by the basic theory of C*-modules
we have  $z = \lim_t e_t z \in Z^2 A \subset Z^2$. 
Thus $Z \cap A \subset Z \cap Z^2$, and so these two sets
are equal.

(4) \  If $c \in Z \cap Z^2$ and
$z \in Z$, then $z c \in Z^2 \cap Z^3 \subset Z^2 \cap Z = J(Z)$;
 and similarly  $ac \in Z \cap Z^2$ for $a \in A$.

(5) \  $J(Z)$ is a closed $A$-$A$-submodule of the
C*-bimodule $Z$.
However the complete $M$-ideals in a
C*-bimodule are exactly the closed $A$-$A$-submodules.
This may be seen for example from the facts in \cite{BEZ}
that the complete $M$-ideals are the left $M$-ideals
which are also right $M$-ideals, and that the
one-sided $M$-ideals in a
C*-bimodule are exactly the closed
one-sided submodules.
\end{proof}

We now return to the center ${\mathcal Z}(E)$
of a $*$-TRO $E$ in a C*-algebra $B$.  

\begin{lemma}  \label{hal}  Let $E$ be a $*$-TRO in a C*-algebra $B$. Then:
 \begin{itemize}
\item [(1)]  ${\mathcal Z}(E)$ is a $*$-TRO in $B$. 
\item [(2)]  If $v \in {\mathcal Z}(E)$
then $v x = x v$ for all $x \in E$.
\item [(3)]  If $v, w \in {\mathcal Z}(E)$
then $v w$ is in the center $C$ of $E^2$.
Also $C {\mathcal Z}(E) \subset {\mathcal Z}(E)$ and 
${\mathcal Z}(E) C \subset {\mathcal Z}(E)$.
\item [(4)]   ${\mathcal Z}(E)$ is a commutative
involutive C*-bimodule over $C$.
\end{itemize}
\end{lemma}
 
\begin{proof}  (1) \  This is immediate. 

 (2) \ Let $u, v, w \in {\mathcal Z}(E)$ and $x \in E$.  Then
$$u v w x = u (w x) v =
(w x) (u v) = (x u) w v = x (w v) u = x w v u.$$
In particular, $(u v)^2 = (v u)^2$.  Setting $v = u^*$,
 and using the unicity of square roots, gives $u u^* = u^* u$.
By the polarization identity we obtain that
 $u v = v u$  for all $u, v \in {\mathcal Z}(E)$.   
Using this together with the last centered equation we obtain 
 $$u v w x = x w v u = x u v w. $$
Since ${\mathcal Z}(E)^3 = {\mathcal Z}(E)$, we obtain (2).

Items (3) and (4) are now easy to check directly (see also \cite{Skei}).
\end{proof}

\begin{lemma} \label{wwmu}  Let $Z$ be a $*$-WTRO in a 
W*-algebra $B$.  Let $J$ be the weak* closure of $Z^2$ in $B$.
We have:
\begin{itemize} \item [(1)]  $Z + J$ 
equals the weak* closure of $Z + Z^2$ in $B$.
\item [(2)]  The multiplier algebra $M(Z^2)$ is $*$-isomorphic to $J$.
 \item [(3)] $M(Z + Z^2) \cong Z + J$ $*$-isomorphically.
\item [(4)]  $J(Z)$ is a W*-subalgebra of $J$, and in particular is
weak* closed in $Z$ too.
\item [(5)]  $Z \cong J(Z) \oplus^\infty J(Z)^\perp$,
where $J(Z)^\perp = \{ z \in Z : z J(Z) = 0 \}$.    \end{itemize}
\end{lemma} 

\begin{proof} 
 Clearly $Z$ and $J$ are 
both contained in the weak* closure of $Z + Z^2$,
and hence so is their sum.  Note that $Z$ is an involutive
$J$-subbimodule of $B$.  Hence by Corollary \ref{iscl} (2),
$Z + J$ is weak* closed.  Hence $Z + J$ contains
the weak* closure of $Z + Z^2$.  This yields (1).
 
(2) \  (See Appendix in \cite{OIAN}). 
By Corollary \ref{iscl} (2),
 $Z + J$ is a W*-algebra,
and it may therefore be represented 
faithfully and normally as a von Neumann algebra in $B(H)$ say. 
 We claim that in this
representation, $Z^2$ acts nondegenerately
on $H$.  For if $\eta$ is a unit vector in the orthocomplement
of $Z^2 H$, and if $\psi = \langle \cdot \eta , \eta \rangle$,
then $\psi$ is a normal state on $Z + J$.  Now $Z^2$ and 
$Z + Z^2$ have a common increasing approximate identity
$(e_t)$, and $e_t \to I_H$ weak*.  Hence 
$\langle \eta , \eta \rangle = \lim \langle e_t \eta , \eta \rangle = 0$,
and so $\eta = 0$.   Thus indeed $Z^2$, and also $Z + Z^2$,
act nondegenerately on $H$.   
Thus we may identify $M(Z^2)$ and $M(Z + Z^2)$ with 
$*$-subalgebras of $B(H)$.  Indeed we may identify them with
$*$-subalgebras of the double commutants of $Z^2$ and 
$Z + Z^2$ respectively in $B(H)$.  Thus, by the
double commutant theorem and (1), we may identify them with
$*$-subalgebras of $J$ and $Z + J$ respectively.   

Conversely, by routine weak* density arguments
$J Z \subset Z$  since $Z$ is weak* closed.  Thus
$J Z^2 \subset Z^2$, and similarly $Z^2 J \subset Z^2$.
Thus $J \subset M(Z^2)$, and so $J = M(Z^2)$.
Similarly $(Z + J)(Z + Z^2) \subset Z + Z^2$
and $(Z + Z^2)(Z + J) \subset Z + Z^2$, so that
$Z + J = M(Z + Z^2)$.  This proves (2) and (3).

For (4), simply note that $J(Z) = Z \cap J$ by Corollary \ref{iscl} (3).

Finally, let $e$ be the identity of the W*-algebra $J(Z)$, and 
$1$ the identity of $J = M(Z^2)$.  If $z \in Z$ then $z = ez + (1-e) z$.
We have $e z \in J(Z)$ (since the latter is an ideal), and 
$(1-e) z \in J(Z)^\perp$.   This proves (5).
\end{proof}

\section{Involutive ternary systems}

We will use the term {\em involutive ternary system}
for the  {\em abstract} version of a  $*$-TRO.
Namely, an involutive ternary system
is a ternary system $X$ possessing an involution 
$*$, 
such that $X$ is completely isometrically isomorphic 
via a selfadjoint (i.e.\ `$*$-linear') ternary morphism,
to a $*$-TRO.   We will give a useful
characterization of these spaces in the next Lemma.
The appropriate morphisms between
involutive ternary systems 
are of course the {\em ternary $*$-morphisms},
namely the  selfadjoint ternary morphisms.

\begin{lemma}\label{*-TRO_char}
Let $X$ be a ternary system possessing an involution
$*$.  Then:   
\begin{itemize} \item [(1)]  (Cf.  \cite[Prop. 2.3]{IKR}) \
$X$ is an involutive ternary system
if and only if 
$[x,y,z]^*=[z^*,y^*,x^*]$ for all $x,y,z\in X$,
and if and only if $\Vert [x_{ij}^*] \Vert = 
\Vert [x_{ji}] \Vert$ for all $n \in \Ndb$ and
$[x_{ij}] \in M_n(X)$.
  \item [(2)]  If $X$ is
an involutive ternary system possessing
 a predual Banach space then
$X$ is isomorphic to a $*$-WTRO via
a weak* homeomorphic completely isometric ternary $*$-isomorphism.
\end{itemize}
\end{lemma}
\begin{proof} (1) \  The one direction of the first
equivalence is obvious.  For the other,
fix a completely isometric ternary morphism $\Phi_0 : X\to B(H)$.
If the hypothesized identity holds, it
is straightforward to check that
$$
\Phi(x):=\begin{pmatrix}0 & \Phi_0(x)\\
                         \Phi_0(x^*)^* & 0
           \end{pmatrix}
$$
is a completely isometric ternary $*$-morphism from 
$X$ onto a $*$-TRO inside $M_2(B(H))$. 

For the second equivalence, again the one direction is
obvious.  For the other, suppose that $X$ is a TRO in $B(H)$,
and that $\tau : X \to X$ is a conjugate linear map satisfying
$\Vert [\tau(x_{ij})] \Vert_n = \Vert [x_{ji}] \Vert$ for 
all matrices.  The map $\theta(x) = \tau(x)^*$ is then a linear
complete isometry from $X$ onto $\{ x^* \in B(H) : x \in X \}$.
Thus by Proposition \ref{prtm} (3), $\theta$ is a ternary 
morphism.   Thus $\tau([x,y,z]) = [\tau(z),\tau(y),\tau(x)]$ for all $x, y, z$.
  
(2) \ If further $X$ has a predual, then by
Zettl's characterization of 
dual TRO's (see \cite[section 4]{ACOT} and
\cite{OIAN}), there is a weak* homeomorphic completely isometric ternary 
morphism $\Phi_0$ from $X$ onto a weak* closed TRO inside $B(H)$. 
 Defining $\Phi$ as in the proof 
of (1) yields the desired 
isomorphism.     \end{proof}

In view of (2) of the Lemma, we define a {\em dual involutive ternary system} 
to be an involutive ternary system with a predual Banach space.
By the Lemma, this is the abstract version of a $*$-WTRO.


Lemma \ref{*-TRO_char} (1) immediately gives another proof of a
result we mentioned after Definition \ref{incb}:

\begin{corollary}  \label{wcor}   Any involutive $C^*$-bimodule
(in the sense of \ref{incb}) is an involutive ternary system.
\end{corollary}

Next we discuss the second dual of a $*$-TRO.                  
If $X$ is an involutive ternary system, then so is $X''$ in a canonical
way. That is, there exists one and 
only one way to extend the involution on $X$
to a weak* continuous involution on $X''$, and with this extended involution
$X''$ is an involutive ternary system. The `existence' here is easy, if $X$
is represented as a $*$-TRO in a C*-algebra $B$, then
$X^{\perp \perp}$ is easily seen by standard arguments
to be a $*$-WTRO in $B''$.  The `uniqueness'
follows by routine weak* density considerations.

\begin{proposition}\label{jinjs1}
Let $Z$ be a $*$-TRO in a C*-algebra $A$, and set  $E = Z''$, which we
may identify with a subspace of the W*-algebra $A''$.
Then we have:
\begin{itemize}
\item [(1)]
$(Z^2)^{\perp \perp}$ in $A''$ equals  the weak* closure $N$ of $E^2$; and
may also be identified with $M(E^2)$, the multiplier algebra of $E^2$.
\item [(2)]
$E + N$ is weak* closed in $A''$.
 \item [(3)]
$J(Z)^{\perp \perp} = J(E)$. Thus $J(Z)$ is weak* dense in $J(E)$.
\item [(4)]  $J(Z)_+$ is weak* dense in $J(E)_+$.
  \end{itemize}
\end{proposition}
 
\begin{proof}  (1) \ Let $J = (Z^2)^{\perp \perp}$.  
Clearly $J \subset N$.  Conversely, by routine
weak* density considerations, $Z E \subset J$,
and so $E^2 \subset J$, giving $N = J$.  Also by routine
weak* density considerations. Since $E$ is an involutive ternary 
$J$-subbimodule of $A''$, $J = M(E^2)$ by 
Lemma \ref{wwmu} (2).

(2) \   Follows from  Corollary \ref{iscl} (2)
(or from  \cite{HWW} I.1.14).

(3)  \ Follows from the above and \cite{HWW} I.1.14.

(4) \ Follows from (3) and the Kaplansky density theorem.
\end{proof}
 
In the following result we use 
complete $M$-ideals, in the sense of
Effros and Ruan \cite{Mide}.
We omit most of the 
proof, which may be found explicitly in \cite{BaTim}, and is based in
part on earlier results by Horn and Kaup
\cite{Horn2,Kaup2,Kaup}.  Of course in our particular setting,
the proofs simplify considerably, but we have resisted the
temptation to include such.

\begin{lemma} \label{3mi}  Let $Z$ be a $*$-TRO.
\begin{enumerate}
\item [(1)]  If $Z$ has a predual then the 
selfadjoint complete $M$-summands in $Z$ are exactly the 
weak* closed ternary $*$-ideals, and these
are exactly the subspaces $p Z$ for a
(necessarily unique) central projection $p$ in the
multiplier algebra $M(Z^2)$ such that
$p z = z p$ for all $z \in Z$.   
  \item [(2)]   The
selfadjoint complete $M$-ideals in $Z$ are exactly the 
ternary $*$-ideals, and these
are exactly the subspaces of the form 
$Z_p = \{ z \in Z : p z = z \}$ for a
projection in the second dual of
$Z^2$ such that $p z = z p$ for all $z \in Z$.  
In fact $p$ may be chosen to have the additional property
of being an open projection (see \cite[3.11.10]{Ped} for the 
definition of this) 
in the center of the second dual of 
$Z^2$, and with this qualification
the correspondence $p \mapsto Z_p$ is bijective.
\end{enumerate} 
Indeed in (1) and (2) above the word `complete' may be dropped.
\end{lemma}

\begin{proof}     
(1) \  If $Z$ has a predual
then we may assume that $Z$ is a $*$-WTRO in a
von Neumann algebra $B$, by Lemma 
\ref{*-TRO_char} (2). 
If $E$ is a weak* closed ternary $*$-ideal in  $Z$, then by the 
facts mentioned in the proof of Corollary \ref{iscl} (5),
$E$ is a complete $M$-ideal in $Z$.  Since 
$E$ is weak* closed, $E$ is a complete $M$-summand
by the basic theory of $M$-ideals \cite{HWW}.

(2) \  Again this is nearly all in the cited literature.  If
$E$ is a ternary $*$-ideal, then by routine arguments, 
$E^{\perp \perp}$ is  a 
weak* closed ternary $*$-ideal of $Z''$.
By (1) it is a selfadjoint 
complete $M$-summand, of $Z''$.  Thus $E$ is a selfadjoint complete $M$-ideal.
For the uniqueness,
 suppose that $p, r$ are two open projections
with the property that if $z \in Z$ then
$p z = z$ if and only if $r z = z$.
If $a \in Z^2$  and if $p a = a$ then
$p a z = a z$ for all $z \in Z$.  Thus $r a z = a z$ for all $z \in Z$,
so that $r a = a$.  Let $a_t $ be an increasing net in $Z^2$ 
converging in the weak* topology to $p$.  Then 
$p a_t = a_t$, so that $r a_t = a_t$.  In the limit
$r p = p$.   By a similar argument $p r = r$.  So $r = p$.     
\end{proof}

\section{Orderings on ternary systems} \label{order} 

In this section,
 we develop the basic theory of ordered ternary systems.
  
\begin{lemma}  \label{jzp}  If $Z$ is a $*$-TRO, then  
$J(Z) = $ {\rm Span}$(Z_+)$, and $J(Z)_+ = Z_+$.
Also, $J(M_n(Z)) = M_n(J(Z))$, and $M_n(Z)_+ =
M_n(J(Z))_+$, for every $n \in \Ndb$.
\end{lemma} 

\begin{proof}    Clearly $J(Z)_+ \subset Z_+$.
  If $z \in Z_+$ then since $z^2 \in Z^2$ we must have
$z \in Z^2$ (square roots
in a C*-algebra remain in the C*-algebra).
Thus  $Z_+ \subset Z \cap Z^2 = J(Z)$, and so $J(Z)_+ = Z_+$.
That $J(Z) = $ {\rm Span}$(Z_+)$
 follows since $J(Z)$ is spanned by $J(Z)_+ = Z_+$. Finally,
$$J(M_n(Z)) = M_n(Z) \cap M_n(Z)^2 =
M_n(Z) \cap M_n(Z^2) = M_n(J(Z)),$$
 and so $M_n(Z)_+ = J(M_n(Z))_+ = M_n(J(Z))_+$.
   \end{proof}      

The following is 
clear from Lemmas \ref{wwmu} (5)
and \ref{jzp}:

\begin{corollary} \label{spl}  If $E$ is a $*$-WTRO
then $E$ is ternary $*$-isomorphic, via a complete order isomorphism,
to $M \oplus^\infty Z$, where $M$ is a W*-algebra, and
$Z$ is a trivially ordered $*$-WTRO.
\end{corollary}

\begin{corollary} \label{contof}  Let
$Z$ be a $*$-TRO, and let
$\theta : Z \to B$ be a positive ternary $*$-morphism
into a C*-algebra $B$.  Then
\begin{itemize} \item [1)]  $\theta$ is completely positive.
 \item [2)] $\theta$ restricts to a
$*$-homomorphism from $J(Z)$ into $B$.
Indeed if $\pi : Z^2 \to B$ is the $*$-homomorphism
canonically associated with $\theta$ (see Proposition 
\ref{prtm} (6)), then
$\theta = \pi$ on $J(Z)$.
\item [3)]  If $A = Z + Z^2$, then
$\theta$ is the restriction
to $Z$ of a $*$-homomorphism from $A$ into $B$.
Hence $\theta$ may be extended further to a
unital $*$-homomorphism from a unitization
of $A$ into a  unitization of $B$.
\item [4)]  If $\theta$ is an order embedding then
it is  a complete order embedding.  In this case,
if $W$ is the $*$-TRO $\theta(Z)$,
then $\theta$
restricts to a $*$-isomorphism of $J(Z)$ onto $J(W)$,
and $\theta$ is the restriction of 
a $*$-isomorphism between $Z + Z^2$ and $W + W^2$.
\end{itemize}
\end{corollary}
 
\begin{proof}  2) \  Consider
the restriction of $\theta$ to $J(Z)$.
Claim: a positive ternary $*$-morphism $\psi$ between C*-algebras is
 a $*$-homomorphism.  To see this note that
by going to the second dual we may assume that the
C*-algebras are unital.  Then if $\psi(1) = v$ it is
clear that $v$ is a positive partial isometry.  
That is, $v$ is a projection.  
Also $v \psi(\cdot)$ is a $*$-homomorphism,
and $\psi = v^2 \psi(\cdot)$
in this case.  Since $v^2 = v$ this proves the claim.
Using the claim, for $a, b \in J(Z)$
we have $\pi(a^* b) = \theta(a)^* \theta(b) = \theta(a^*b)$.
Thus $\theta = \pi$ on $J(Z)$.
 
3) \ Let $W = \theta(Z)$.
Define $\tilde{\theta} :
Z + Z^2 \to W + W^2$ by $\tilde{\theta}(x + yz) = \theta(x) +
\theta(y) \theta(z)$,
if $x, y, z \in Z$.  Namely $\tilde{\theta}$ is the unique
linear extension of both $\theta$ and the $*$-homomorphism
$Z^2 \to W^2$ associated with $\theta$.   It is easy to check
using 2) that
$\tilde{\theta}$ is a well defined $*$-homomorphism on $Z + Z^2$.
Since $Z + Z^2$ is closed,  $\tilde{\theta}$ is contractive
and will extend to a unital
 $*$-homomorphism between the unitizations.
 
1) \ Follows from 3).
 
4)  \ Since $\theta$ is an order embedding it maps $Z_+$ onto $W_+$.
By Lemma \ref{jzp} 
we see that $\theta$  maps $J(Z) = $ Span $(Z_+)$ 
onto $J(W) = $ Span$(W_+)$.  By 2), $\theta$ is
a $*$-homomorphism.  To see the isomorphism
between $Z + Z^2$ and $W + W^2$, use the proof of  
3), and the same construction applied to $\theta^{-1}$.  
\end{proof}

\begin{definition} \label{odef}
An {\em ordered ternary system} is an involutive  ternary system $X$ 
which is an  ordered operator space too; the cone 
on $X$ will be called an {\em ordered operator space cone}.
On the other hand we shall use the term {\em naturally ordered ternary 
system} for an ordered ternary system such that $X$ is completely 
order isomorphic via a ternary $*$-isomorphism to a $*$-TRO
$Y$, where $Y$ is given its relative cones inherited from 
its containing C*-algebra. 
The associated ordering (resp.\ positive cone in $X$,
positive cones in $M_n(X)$) will be referred to 
as a {\em natural ordering} (resp. {\em natural cone},
{\em natural matrix cones}) on $X$.
Similarly, we refer to  {\em natural dual orderings} 
and  {\em natural dual cones} on a dual ternary system;
this corresponds to the ordering or cone pulled back from
a $*$-WTRO via a weak* homeomorphic ternary $*$-isomorphism.
\end{definition}
 
Note that for every involutive ternary system,
the trivial ordering is a natural ordering.
Indeed given a $*$-TRO $Y$ in  a C*-algebra $B$,
one may replace $Y$ with the isomorphic subspace of $M_2(B)$ consisting 
of matrices of the form 
$$
\begin{pmatrix}
0 & y \\
y & 0 
\end{pmatrix} , \qquad y \in Y. 
$$

If $X$ and $Y$ are $*$-TRO's (resp.\ $*$-WTRO's)
in two C*-algebras (resp.\ W*-algebras) $A$ and $B$,
then $X \oplus^\infty Y$ is a $*$-TRO (resp.\ $*$-WTRO) in 
$A \oplus^\infty B$.  It is easy to see that
${\mathcal Z}(X \oplus^\infty Y) = {\mathcal Z}(X)
\oplus^\infty {\mathcal Z}(Y)$.  More generally,

\begin{proposition} \label{dirs}  The `$L^\infty$-direct 
sum' $\oplus^\infty_i \; X_i$ of a family of involutive ternary
systems (resp. naturally ordered ternary
systems, dual involutive ternary
systems, naturally dual ordered ternary
systems) is again an involutive ternary
system (resp. a  naturally ordered ternary
systems, dual involutive ternary system,
naturally dual ordered ternary
system).  Moreover, ${\mathcal Z}(\oplus^\infty_i \; X_i)
= \oplus^\infty_i \; {\mathcal Z}(X_i)$.
\end{proposition}

This is the usual `$L^\infty$'-direct sum, with the
obvious involution (resp. and positive cones). 

The following definition is a basic JBW*-triple construct
(connected with the important notion of `Pierce decompositions',
etc.):
 
\begin{definition} \label{defju}
Let $u$ be a selfadjoint tripotent (that is,
$u = u^* = u u u$) 
in ${\mathcal Z}(E)$, where $E$ is
an involutive ternary system.  We define $J(u) = u E u$,
also known as the Pierce $2$-space.
By the {\em canonical product} on $J(u)$ 
we mean the product $x \cdot y = x u y$ for $x, y \in J(u)$.  It is well
known, and easy to check,
that with this product, and  with the usual involution,
$J(u)$ is a C*-algebra with identity $u$.
The positive cone in this C*-algebra will be
written as ${\mathfrak c}_u$.
\end{definition} 

It is also easy to see that $J(u)$ is a ternary $*$-ideal
in $E$; and that if $E$ is a $*$-WTRO then $J(u)$ is also  weak* closed.
The following gives some alternative
descriptions  of ${\mathfrak c}_u$:
 
\begin{lemma} \label{spl2}
If $E$ is a $*$-TRO, and if $u$ is a
selfadjoint tripotent in ${\mathcal Z}(E)$,
then
$${\mathfrak c}_u = 
\{ e u e^* : e \in E \}  = \{ x \in E : u x \in (E^2)_+ , x = u x  u \} 
.$$
\end{lemma}
 
\begin{proof}  Clearly ${\mathfrak c}_u =
\{ e u e^* : e \in E \}$,  
which in turn is contained in the right hand set.
Conversely, if $u x \in (E^2)_+$ then there is a net whose terms
are of the form $\sum_k e_k e_k^*$, converging to $ux$.
Multiplying the net by $u$, we see that if $u^2 x = x$ then
$x$ is in the closed convex hull of $\{ e u e^* : e \in E \}$.
However the latter set is the positive cone of $J(u)$
with its canonical C*-algebra structure, and hence is closed and convex.
\end{proof}
 
Before characterizing  natural cones on $*$-TRO's,
we will have to tackle the $*$-WTRO case.  
Much of this is in some sense a rephrasing of JBW*-triple facts and
techniques (a `real form' version of facts in e.g.\
\cite[Section 3]{Bat}).   Since we could not find the version we need 
in the literature, and since it will not take us very long,
we include the quick proofs.     
       
 \begin{lemma}\label{w*-ideal_char} Suppose that $J$ is a
weak* closed ternary $*$-ideal $J$ of a $*$-WTRO $E$.
Suppose further that $J$ is ternary $*$-isomorphic to a W*-algebra 
$\fM$ via a mapping $\Phi:\fM\to J$. If $1$ denotes the unit of $\fM$, then 
$u=\Phi(1)$ is a selfadjoint tripotent in ${\mathcal Z}(E)$, and $J=J(u)$.

Conversely, if $u$ is a selfadjoint tripotent in ${\mathcal Z}(E)$,
then the  ternary $*$-ideal $J(u)$ 
of $E$ is ternary $*$-isomorphic to a W*-algebra, via an isomorphism 
taking $u$ to $1$.
\end{lemma}

\begin{proof} 
Clearly, $u=\Phi(1)$ is a selfadjoint tripotent.
Since $\Phi$ is a  ternary $*$-morphism, we have
 $z = u uz = uzu  = zuu$ for all $z \in J$. It follows that $J\sub uEu$.
On the other hand, $uEu = u^3 E u \sub u^2 E \sub J$, so that $J = u E u = u u E$.
Similarly $J = E u u$.   Thus the maps
$e\mapsto uue$ and $e\mapsto ueu$ are norm-one projections from $E$ 
onto $J$. But $J$ is 
a complete $M$-summand by Lemma \ref{3mi} (1),
and hence an $M$-summand,
 and the projections onto such subspaces are unique (see 
e.g.\ \cite[Proposition I.1.2]{HWW}). Hence, $ueu=uue=euu$ for all $e\in E$. 
Thus
$$
ue = u^3 e = u^2 e u = e u^3 = e u,$$ 
and hence $u \in {\mathcal Z}(E)$. 

It is, on the other hand, well-known (and in any case an
easy exercise) to check
that any weak* closed subspace of the form 
$J(u)$ is ternary $*$-isomorphic to the W*-algebra
$u E$.  The induced product on $J(u)$ is just its
canonical product mentioned above
the lemma.
\end{proof}

We now link natural dual orderings with `central tripotents':

\begin{theorem}\label{w*-natcone_char}
Let $E$ be an  involutive ternary system.
A given cone $E_+$ in $E$ (resp.\ $M_n(E)_+ \subset M_n(E)$
for all $n \in \Ndb$) is a natural dual cone
(are natural dual matrix cones) if and only if 
for some selfadjoint tripotent $u\in 
{\mathcal Z}(E)$,
$$
E_+= \{ eue^* : e\in E \}
$$
(resp. \
$$M_n(E)_+ = \{ [\sum_k e_{ik} u e_{jk}^*] : [e_{ij}] \in 
M_n(E) \} , \qquad n \in \Ndb).$$
\end{theorem}

\begin{proof}  ($\Leftarrow$) \
Fix a selfadjoint 
tripotent $u\in {\mathcal Z}(E)$ with $E_+= \{ eue^* : e\in E \}$. 
We may assume that $E$ is a $*$-WTRO in a W*-algebra (but will
not care about the induced order from this W*-algebra).
Note that $e u e^* = e u^5 u^* = (ueu) u (ue^* u)$, since
$u \in {\mathcal Z}(E)$.  It then follows that 
$E_+= {\mathfrak c}_u$, in the notation of Lemma \ref{spl2},
and regarding $J(u)$ as a 
W*-algebra as in the proof of  
Lemma \ref{w*-ideal_char}.   A similar formula clearly holds
for $M_n(E)_+$, for all $n \in \Ndb$.
Fix a $*$-isomorphism $\Phi_1$ of
 $J(u)$ onto a W*-subalgebra of some 
space $B(H_1)$. Then $\Phi_1$ is a ternary 
$*$-isomorphism which is a complete order 
isomorphism onto its image.   Now $J(u)$ is an $M$-summand of
$E$ as we observed in the proof of Lemma \ref{w*-ideal_char}.
Indeed $J(u) = \{ u u x : x \in E \}$.  Let $J(u)^\perp$ be the
`complementary  $M$-summand' $\{ x-uux : x\in E \}$.
This is a ternary 
$*$-ideal of $E$ too, and in fact is the 
orthocomplement of $J(u)$ in the C*-module sense too
\cite{W-O}.   The map $x \mapsto (u^2 x , x-uux)$ is
a completely isometric ternary $*$-isomorphism from
$E$ onto $J(u) \oplus^\infty J(u)^\perp$ (see 
Proposition \ref{dirs}).  It is also easy to check that
this map is a weak* homeomorphism, and is a
complete order isomorphism when $J(u)^\perp$ is endowed 
with the trivial order.   Let $\Phi_2 : J(u)^\perp \to B(H_2)$
be a weak* continuous ternary $*$-isomorphism
onto a trivially ordered $*$-WTRO in $B(H_2)$ (see the proof of 
Lemma \ref{*-TRO_char}).
Then $\Phi=\Phi_1\oplus\Phi_2$ is a weak* continuous
ternary $*$-isomorphism from $E$ onto a $*$-WTRO, which also is a complete 
order isomorphism onto its image.

($\Rightarrow$) \ Suppose that $E$ is a $*$-WTRO in $B(H)$.  By 
Lemma \ref{jzp}, the cone of $E$ is the cone of the W*-algebra $J(E)$. 
By Corollary \ref{iscl} (5) and Lemma \ref{wwmu} (4),
$J(E)$ is a weak* closed ternary $*$-ideal of $E$. 
So if $u$ is the unit of $J(E)$ then, by 
Lemma~\ref{w*-ideal_char}, $u\in {\mathcal Z}(E)$, and 
$J(E) = J(u)$.
If $e \in E $ then  $eue^*=(euu)u(uue^*)$; and 
also $euu,uue^*\in J(E)$ since the latter is a 
ternary $*$-ideal. Hence
$$
E_+=J(E)_+= \{ eue^* : e\in J(E) \}=\{ eue^* : e\in E \}.
$$
Similar arguments apply to $M_n(E)_+$, which 
equals $M_n(J(E))_+$ by Lemma \ref{jzp}.
 \end{proof}

An important remark is that
by Lemma \ref{jzp} and Corollary
\ref{contof} (4) for example, a natural cone on an involutive
ternary system $Z$ completely determines the associated cone
on $M_n(Z)$.  Nonetheless it is convenient to have the
above explicit description of these cones.

Most of the next two results
are essentially in e.g.\ \cite{ER1,ER4}, \cite[Section 3]{Bat},
\cite{IKR}.
For completeness we include proof sketches:
 
\begin{lemma}\label{z(e)w*}
Let $E$ be a $*$-WTRO. Then
\begin{itemize}
\item[(1)]
The $*$-WTRO ${\mathcal Z}(E)$
possesses a selfadjoint unitary $u$,
such that ${\mathcal Z}(E)\sub J(u)$.  Moreover, equipping
$J(u)$ with its canonical product,
${\mathcal Z}(E)$ is a commutative unital W*-subalgebra
of $J(u)$.
\item[(2)]
An element $v\in {\mathcal Z}(E)$ is a selfadjoint tripotent 
if and only if it is such in $J(u)$, where $u$ is as in {\rm (1)}.
\item[(3)]  If $v$ is a selfadjoint tripotent in ${\mathcal Z}(E)$,
 then there exists a  selfadjoint unitary $w$ in ${\mathcal Z}(E)$
such that $v w v = v$.
\end{itemize}
\end{lemma}

\begin{proof}
(1) \ 
We observed earlier that ${\mathcal Z}(E)$ is a $*$-TRO in $E$, and 
indeed may be viewed as a $*$-WTRO in a W*-algebra $M$.
The unit ball of ${\mathcal Z}(E)_{sa}$ is
weak* closed in the weak* closed subspace $M_{sa}$ of $M$.
By the Krein-Milman theorem,
the unit ball of ${\mathcal Z}(E)_{sa}$ contains an extreme point $u$.   
A standard  Urysohn lemma argument shows 
that $u^2 = 1_A$ in $A$ (one shows that every character $\chi$ on $A$ 
satisfies $\chi(u^2) = 1$). 
Thus if $x \in {\mathcal Z}(E)$ then $x = u u x \in J(u)$,
and so ${\mathcal Z}(E)\sub J(u) $ as desired.   
Since ${\mathcal Z}(E)$ is a $*$-TRO it is clear that
${\mathcal Z}(E)$ is closed under the canonical product of
$J(u)$.  Now it is clear that 
${\mathcal Z}(E)$ is a commutative W*-subalgebra
of $J(u)$, and the unit $u$ lies in this
subalgebra.

(2) \ This is clear from the fact that $u$ is unitary, since
$vuvuv= u^2 vvv = vvv$ if $v \in {\mathcal Z}(E)$.

(3) \  Let $v$ be a selfadjoint tripotent in ${\mathcal Z}(E)$,
and let $u$ be as in
(1).  Then $w = v + u - v u v$ is the
desired unitary.     
\end{proof}

\begin{proposition}\label{bij_order_corr}
Let $E$ be a dual involutive ternary system.
\begin{itemize}
\item[(1)]
There is a bijective correspondence between natural dual
orderings on $E$ and 
selfadjoint tripotents of ${\mathcal Z}(E)$. This correspondence respects the 
ordering of cones by inclusion on one hand, and the ordering of 
selfadjoint tripotents given by $u_1\leq u_2$ if and only if $u_1u_2u_1=u_1$. 
\item[(2)]
If $v\in {\mathcal Z}(E)$ is a selfadjoint
tripotent, then the following conditions are equivalent:
\begin{itemize}
\item[(a)]
$v$ is an extreme point of ${\mathcal Z}(E)_{sa}$.
\item[(b)]
$v$ is a selfadjoint unitary of ${\mathcal Z}(E)$.
\item[(c)]
$v$ is maximal with respect to the order defined in (1).
\item[(d)]
${\mathcal Z}(E)\sub J(v)$.
\item[(e)]  $J(v)$ contains every selfadjoint unitary of ${\mathcal Z}(E)$.
\end{itemize}
\item[(3)]    The maximal natural dual orderings
on $E$ are in 1-1 correspondence with the
selfadjoint unitaries of ${\mathcal Z}(E)$.
\item[(4)] Any natural dual ordering on $E$ is majorized by a maximal 
natural dual ordering.
\end{itemize}
\end{proposition}
\begin{proof}  (1) \  It is an easy exercise (or see the above cited 
papers) that $\leq$ is a partial ordering.
Write ${\mathfrak c}_1$ and ${\mathfrak c}_2$ for the
corresponding natural dual cones (see Theorem \ref{w*-natcone_char}).
  If ${\mathfrak c}_1 \subset {\mathfrak c}_2$,
then by Lemma \ref{spl2} we have  
 $u_1 u_2 \geq 0$, and $u_1^2 u_2^2 = u_1^2 = u_1^4$.
Putting these together we have
$u_1 u_2 = u_1^2$, so that $u_1 \leq u_2$.  Conversely, 
if $u_1 \leq u_2$ then $e u_1 e^*  = e u_1 u_2 u_1 e^* \in {\mathfrak c}_2$.
Hence ${\mathfrak c}_1 \subset {\mathfrak c}_2$.
Now (1)  follows readily from Theorem~\ref{w*-natcone_char}.

(2) \  These are well-known; but also quite clear: 
For the equivalence of (a) and (b) 
see e.g.\ the hints in the proof of Lemma \ref{z(e)w*} (1).
To show that (b) implies (c) suppose that $v'\in {\mathcal Z}(E)$ is a selfadjoint 
tripotent with $vv'v=v$. By assumption, $vv'v=v'$. So $v = v'$.
That (c) implies (b) follows from Lemma \ref{z(e)w*} (3).
Clearly (d)  implies (e); and the converse is routine.

(3)  \ Follows from (1) and (2).   

(4) \  Follows from (1), (3), and Lemma \ref{z(e)w*} (3).
\end{proof}

It is interesting to deduce structural
facts about orderings on a $*$-WTRO $E$ from the algebraic structure of the
partially ordered set ${\mathcal S}(E)$ of selfadjoint tripotents in $Z(E)$
(see e.g.\ \cite{AkeP,Bat,ER1}).
In fact ${\mathcal S}(E)$ is
not a lattice: the `sup' in general does not make sense without further
hypotheses.  However `infs' are
beautifully behaved.  Indeed in our situation, it is rather clear that
for any family of natural dual cones $({\mathcal C}_i)$ on $E$,
there is a natural dual cone which is the
largest natural dual cone on $E$ contained in every ${\mathcal C}_i$.
Indeed this cone is just $\cap_i {\mathcal C}_i$.
  To see that this indeed is a  natural dual cone, we fix for each
$i$ a ternary order isomorphism $T_i$ from $(E,{\mathcal C}_i)$ onto a
$*$-WTRO $E_i$.  Then define  $T(x) = (T_i(x))_i$, this is
a ternary order isomorphism from
$(E,\cap_i {\mathcal C}_i)$ into $\oplus^\infty_i E_i$.
By Proposition \ref{bij_order_corr} it follows that for every
family of selfadjoint tripotents in ${\mathcal S}(E)$, there is
a greatest lower bound in  ${\mathcal S}(E)$.  In fact there
is a convenient explicit formula.
If $u, v$ are selfadjoint tripotents in ${\mathcal S}(E)$,
then the greatest lower bound of $\{ u , v \}$ in ${\mathcal S}(E)$ is
\begin{equation} \label{meet}
u \wedge v = \frac{1}{2}(u v u + v u v)
\end{equation}
We leave the  proof of this as a simple algebraic exercise 
(see also \cite{AkeP}).

 Continuing the discussion in the last 
paragraph, clearly $u \in
{\mathcal S}(E)$ if and only if $-u \in {\mathcal S}(E)$;
however the reader should be warned that 
 the natural ordering $\leq$ on tripotents
(see (1) of the previous result) behaves a little unexpectedly:
$u \leq v$ if  and only if
$-u \leq -v$.

It is worth stating separately the fact
that if $E$ is a $*$-WTRO, and if $u$ is the selfadjoint
tripotent in ${\mathcal Z}(E)$ corresponding (as in the last
proposition) to
the given order on $E$, then $J(E)$, and the usual product on
$J(E)$ (recall that we pointed out in Section \ref{ICB}
that this is a C*-algebra), may be recaptured in terms
of the ternary structure and the tripotent $u$, as
precisely the
canonical product on $J(u)$.  This observation is rather
trivial (see the proof of Theorem 
\ref{w*-natcone_char}), but is useful
when $E$ is not given as a $*$-WTRO, but instead as an
abstract involutive ternary system.

\begin{corollary} \label{maxcnd} The given cone on a $*$-WTRO $E$ is maximal among 
all the natural dual cones on $E$, if and only if $E$ is ternary order 
isomorphic to $M \oplus^\infty F$, where 
$M$ is a W*-algebra and $F$ is a $*$-WTRO for which
${\mathcal Z}(F) = \{ 0 \}$.
\end{corollary} \begin{proof} 
$(\Rightarrow)$ \   By Corollary \ref{spl},
 $E \cong J(E) \oplus^\infty 
J(E)^\perp$.
Take $u$ to be the associated  selfadjoint  unitary  
in ${\mathcal Z}(E)$ (see Proposition \ref{bij_order_corr} (4)).
Note that ${\mathcal Z}(J(u)^\perp) \subset {\mathcal Z}(E) \cap J(u)^\perp
= (0)$, by 2 (d) in  Proposition \ref{bij_order_corr}. 

$(\Leftarrow)$ \  In this case 
${\mathcal Z}(E) = {\mathcal Z}(M) \subset J(E) = M$.
Note that $(1,0)$ is a selfadjoint central unitary with $J(u) = M$.   
Now appeal to Proposition \ref{bij_order_corr}.
  \end{proof} 

{\bf Remarks.}  1) \  The last result shows that dual involutive
ternary systems have a quite simple canonical form.
Such canonical forms are quite common in the 
JBW*-triple literature.  As we discuss elsewhere however,
$*$-TRO's $F$ with ${\mathcal Z}(F) = \{ 0 \}$ are not necessarily
uncomplicated.  

2) \ 
We shall see later in Proposition \ref{maxdo} that the set of
maximal natural dual cones on $E$ (characterized in the last two results
above),
coincides with the
set of maximal ordered operator space cones on  $E$.

\medskip

We now pass to $*$-TRO's.
We recall from the discussion above Proposition \ref{jinjs1},
that if
$Z$ is a $*$-TRO in a C*-algebra $B$, then $E = Z^{\perp
\perp}$ is  a $*$-WTRO in $B''$.
One may deduce from Proposition  
\ref{jinjs1} (3) that there is a smallest weak* closed cone on $E = Z''$
which contains the cone given on $Z$, namely the weak* closure of this
cone, and this is a natural cone for $Z''$.  We call this the {\em
canonical ordering} or {\em canonical second dual cone}
on $Z''$.  However there may in general
be many other (bigger) natural cones on $Z''$ which induce the same cone
${\mathfrak c}$ on $Z$ (see the examples towards the end of
Section \ref{ilb}).

\begin{proposition}\label{jinjs}
Let $Z$ be a $*$-TRO in a  C*-algebra $A$, and let  $E$
be the $*$-WTRO $Z''$ in the W*-algebra $A''$.
Suppose that $T : Z \to B$ is a $*$-linear map
into a C*-algebra $B$. 
Then:
\begin{itemize}
\item [(1)]
$T$ is completely positive if and only if $T''$ is 
completely positive.
\item [(2)]   If $T$ is a ternary $*$-morphism,
then so is $T''$. 
In this case,
$T$ is a complete order embedding if and only if $T''$ is
a complete order embedding.
\end{itemize}
\end{proposition}
 
\begin{proof} 
(1)  \ By Proposition \ref{jinjs1},
$J(E)_+ = (J(Z)'')_+$.  If $\eta \in (Z'')_+ = J(E)_+ = (J(Z)'')_+$, then by 
a variant of 
Kaplansky's density 
theorem, there exists a net $(x_t)$ in $J(Z)_+$ converging to $\eta$ in 
the weak* topology.  Then $T''(x_t) = T(x_t) \to T''(\eta)$. But $T(x_t)\geq 
0$, and so $T''(\eta) \geq 0$. 
Thus $T''$ is a positive map.   Since $T_n : M_n(Z) \to M_n(B)$ satisfies 
$(T_n)'' = (T'')_n$ we deduce that $T''$ is completely positive. The other
direction follows by simply restricting $T''$ to $Z$.
 
(2) \ The first assertion follows by routine weak* density
arguments.      The second
follows from (1) applied to $T$ and $T^{-1}$.
\end{proof}

\begin{definition} \label{defdc}  Let  $Z$ be an involutive ternary system,
and suppose that $u$ is a selfadjoint tripotent $u$ in ${\mathcal Z}(Z'')$.
If ${\mathfrak c}_u$ is the (natural) cone defined in \ref{defju}, we define 
${\mathfrak d}_u$ be the cone
${\mathfrak c}_u \cap Z$ in $Z$.    We also write
$J_u(Z)$ for the span  of ${\mathfrak d}_u$  in $Z$, and   
${\mathfrak c}'_u$ for the weak* closure of
${\mathfrak d}_u$ in $Z''$. 

We say that a tripotent $u$ is {\em open}
if, when we consider $J(u)$ as a W*-algebra in the canonical
way (see Lemma \ref{w*-ideal_char}),
then $u$ is the weak* limit in $Z''$ of an increasing net
in $J(u)_+ \cap Z$. 
\end{definition}

Later we will give several 
alternative characterizations of open central tripotents, and list some of
their basic properties.  We note that an projection in the second dual
of a C*-algebra is open in the usual 
sense if and only if it is an open tripotent in the
sense above.

\begin{lemma} \label{uxinE}   Let $Z$ be an
involutive ternary system, and let $E = Z''$,
also an involutive ternary system in the canonical way. 
Let $u$ be a selfadjoint tripotent  in ${\mathcal Z}(E)$.
We have:
  \begin{itemize} \item [(1)]  ${\mathfrak c}'_u \subset {\mathfrak c}_u$.
 \item [(2)]  ${\mathfrak d}_u$  is a natural cone on 
$Z$.  That is, there is a surjective  
ternary order $*$-isomorphism $\psi$ from $Z$ with cone ${\mathfrak d}_u$
to $W$,
for a $*$-TRO $W$.  Also $\psi''$ is a 
ternary order $*$-isomorphism from $Z''$ with cone ${\mathfrak c}'_u$,
 onto $W''$ with its canonical second dual cone.
\item [(3)]   ${\mathfrak c}'_u = {\mathfrak c}_v$ for a
selfadjoint tripotent $v \in {\mathcal Z}(E)$ with $v \leq u$.
\item [(4)]  $J_u(Z)$ is
a C*-subalgebra of $J(u)$, the latter
regarded as a C*-algebra in the canonical
way (see Lemma \ref{w*-ideal_char}).
Also, ${\mathfrak d}_u$  is
the positive cone of this C*-algebra $J_u(Z)$.
\item [(5)]  If $Z$ is a $*$-TRO in a C*-algebra $B$ say (we do not care 
about the ordering induced on $Z$ by $B$),
then $u {\mathfrak d}_u \subset {\mathfrak d}_u^2 \subset
J_u(Z)^2 \subset Z^2$.  Thus $u J_u(Z)  \subset J_u(Z)^2$.
\end{itemize}
\end{lemma}  
 
\begin{proof}   (1) \ By Theorem \ref{w*-natcone_char} and
Lemma  \ref{spl2}, $E$ may be regarded as
a $*$-WTRO whose positive cone is ${\mathfrak c}_u$, and this
is weak* closed in $E$.   Then (1) is clear.

(2) \  Continuing with the argument in (1),
regarding  $Z$ as a $*$-TRO inside $E$, ${\mathfrak d}_u$ is
just the inherited cone.   Thus it is a natural cone.
By  Proposition~\ref{jinjs1} (3) and Kaplansky's density
theorem, ${\mathfrak c}'_u$ coincides with the 
canonical second dual cone induced from ${\mathfrak d}_u$.

(3) \  Follows from the last part of (2), and Theorem
\ref{w*-natcone_char} and Proposition \ref{bij_order_corr}.

(4) \ This follows from the above, and Lemma \ref{jzp}.

(5) \  By (4), any element $x$ of ${\mathfrak d}_u$ may be written as $y u y$ for
some $y \in {\mathfrak d}_u$.  
Thus $u x = u^2 y y = y y \in {\mathfrak d}_u^2 \in J_u(Z)^2$. 
Then (5) is clear.
\end{proof}           

\begin{theorem}\label{order_char}
Let $Z$ be an involutive ternary system. Then
\begin{itemize}
\item[(1)]
The natural cones on $Z$ are precisely the
cones
${\mathfrak d}_u  = {\mathfrak c}_u \cap Z$
(notation as in \ref{spl2}), for $u$ a selfadjoint open tripotent
in ${\mathcal Z}(Z'')$.                         
A similar fact holds for the  natural matrix cones on $Z$.
\item[(2)]
There is a bijective correspondence between natural orderings on $Z$, and 
selfadjoint open tripotents in ${\mathcal Z}(Z'')$.
 This correspondence respects 
the ordering of the cones by inclusion on one hand, and the ordering of 
selfadjoint tripotents given by $u_1\leq u_2$ if
and only if $u_1u_2u_1=u_1$. 
 \end{itemize}
\end{theorem}

\begin{proof}
(1) If $Z$ is a $*$-TRO, 
then so is $E = Z''$.
Applying Theorem~\ref{w*-natcone_char} to the canonical
ordering on $E$ yields a 
selfadjoint tripotent $u$ with 
$J(u) = J(E)$, and $E_+ = {\mathfrak c}_u$.
  That $u$, the identity
of the W*-algebra $J(E)$, is open 
follows from Proposition~\ref{jinjs1} (3).  Thus
 $Z_+$ has to be of the announced form. 
A similar argument applies to $M_n(Z)_+$.

The converse follows from Lemma \ref{uxinE} (2).

(2) \ Suppose that $u, v$ are two selfadjoint open tripotents
in ${\mathcal Z}(E'')$.  Set $A = \{eue^* : e\in E'' \} \cap E$,
and $B$ the matching set for $v$.  Suppose first that 
$A \subset B$.
By definition $u$ is a weak* limit of 
a net in $A$, and hence in $B$.  Since ${\mathfrak c}_v$ is
weak* closed (being the cone of a W*-algebra),
we have that $u$ is in ${\mathfrak c}_v$.
It follows as in the proof of Proposition
\ref{bij_order_corr} (1), that $u \leq v$.
Conversely, if $u \leq v$, then the 
second part of Proposition
\ref{bij_order_corr} (1) shows that 
$A \subset B$.
\end{proof}

We now give some equivalent, and often more useful,
characterizations of selfadjoint open central tripotents.  But first we will
need one or two more definitions and facts.
Suppose that  $J$ is a C*-ideal in an involutive ternary system $Z$,
with  $\psi : B \to J$ the surjective triple $*$-isomorphism,
where $B$ is a C*-algebra.   Suppose that $Z$ is
a $*$-TRO in a C*-algebra $A$, and that
$E = Z''$ is regarded as a $*$-WTRO in the W*-algebra $M = A''$.
Then $\psi'' : B'' \to J^{\perp \perp}$ is a
surjective triple $*$-isomorphism onto a
weak* closed triple ideal in $E = Z''$.
Let $\psi''(1) =
u$, by Lemma \ref{w*-ideal_char} this is a selfadjoint
tripotent in ${\mathcal Z}(E)$ and $J^{\perp \perp} = J(u)$.
We call $u$ a {\em support tripotent} for $J$.  If $(e_t)$ is an
increasing approximate identity for $B$, then it is well
known that $e_t \to 1_{B''}$ weak*.  Thus $\psi(e_t) \to u$
weak*.   We deduce that $u$ is an open tripotent in the sense
above.     
Also  $J$ with the product pulled back from $B$ via 
$\psi$, is a C*-subalgebra of $J(u)$, the latter with its canonical
 product (see \ref{spl}).   Clearly $J = J^{\perp \perp}
\cap Z = J(u) \cap Z$.  We will see in Theorem 
 \ref{chopt} that $u^2$ is open in the usual C*-algebraic
sense (from which it
is easy to see that $u^2$ is the projection in Lemma \ref{3mi} (2)).
Finally we remark that the induced natural cone ${\mathfrak d}_u$ 
on $Z$ (see \ref{defdc}), equals $\psi(B_+)$, the 
canonical positive cone for $J$.    To see this note that
$$\psi(B_+) = \psi''(B''_+ \cap B) = \psi''(B''_+) \cap J 
= {\mathfrak c}_u \cap (J(u) \cap Z) = {\mathfrak d}_u.$$

The following theorem is reminiscent of \cite{ER3} Lemma 3.5; but in 
fact only  seems to be formally related to that result.  Indeed 
examples of the type in Section \ref{ilb} show that there 
can exist central selfadjoint tripotents in $Z''$ which are
open in the sense of \cite{ER3,ER2}, but which
are not related to our notion of open tripotent (all selfadjoint
unitaries in the second dual are open in their sense for example).
                                                         
\begin{theorem} \label{chopt}  Let $Z$ be a $*$-TRO, and let $E = Z''$.
Let $u$ be a selfadjoint tripotent in ${\mathcal Z}(E)$.
The following are equivalent: \begin{itemize}
\item [(i)]  $u$ is an open  tripotent.
\item [(ii)]  $u \in {\mathfrak c}'_u$ (notation as in
\ref{defdc}). 
 \item [(iii)]  There is a net $(x_t)$ in $Z$ converging weak* to $u$,
satisfying: $u x_t \geq 0$, $u^2 x_t = x_t$ for all $t$, and 
$(u x_t)$ is an increasing net.
\item [(iv)]  If $J_u(Z)$ is the span of the cone ${\mathfrak d}_u$ (see
\ref{uxinE}) in $Z$, then $\overline{J_u(Z)}^{w*} = J(u)$.
\item [(v)]  $u$ is a support tripotent for a C*-ideal in $Z$.
\item [(vi)]   $u^2$ is an open projection (in the usual sense) in $(Z^2)''$, and 
  $u (J(u) \cap Z) \subset Z^2$. 
 \item [(vii)]  ${\mathfrak c}'_u = {\mathfrak c}_u$ (notation as in \ref{defdc}).
That is, ${\mathfrak c}_u$ is the `canonical second dual cone' induced by
${\mathfrak d}_u$.
  \end{itemize}  \end{theorem}
 
\begin{proof} Clearly (i) is equivalent to (iii), and (iii) implies (ii),
and (vii) implies (ii).  We now show that (ii) implies (vii), (i),  and (iv).  
Using 
Lemma \ref{uxinE}, 
and the first part of the proof
of Theorem \ref{order_char},  we have
${\mathfrak c}'_u = {\mathfrak c}_v$ for a selfadjoint open central tripotent $v$.
Thus if (ii) holds then $u \in {\mathfrak c}_v$, so that 
${\mathfrak c}_u \subset {\mathfrak c}_v
 = {\mathfrak c}'_u$.  But by 
Lemma \ref{uxinE} again, ${\mathfrak c}'_u
\subset {\mathfrak c}_u$.  So ${\mathfrak c}_u = {\mathfrak c}_u'
 = {\mathfrak c}_v$.  Hence $u = v$ is open.  
  By Proposition \ref{jinjs1} (3) and Lemma \ref{jzp},
$\overline{J_u(Z)}^{w*} = J(u)$.
Thus we have verified  (vii), (i),  and (iv).

 (iv) $\Rightarrow $ (i) \  By the proof of Lemma \ref{uxinE} (5),
 we have that $J_u(Z)$ is a C*-subalgebra of $J(u)$.  Then (1) follows
by Kaplansky's density theorem.   

(iv) $\Rightarrow$ (vi) \ Given (iv), 
note that $J(u) \cap Z = J_u(Z)$.  Hence by
Lemma \ref{uxinE} (5),  $u (J(u) \cap Z) = u J_u(Z) \subset  J_u(Z)^2 =
(J(u) \cap Z)^2 \subset Z^2$.
  Since (iv) implies (i), $u^2$ is a weak* limit
of an increasing net in $u(J(u) \cap Z) \subset Z^2$. Thus 
$u^2$ is open.

(v) $\Rightarrow $ (i)  \ We saw this in the discussion above Theorem
 \ref{chopt}.

(vi)  $\Rightarrow $ (v) \ Let $J = J(u) \cap Z$, this is a
ternary $*$-ideal in $Z$, and if (vi) holds then it is
a C*-ideal in $Z$.  Indeed $J$ is a C*-subalgebra of the W*-algebra $J(u)$, 
the latter with its canonical product.  By the discussion  above Theorem
 \ref{chopt}, if $(e_t)$ is an increasing approximate identity
for this C*-subalgebra, then $e_t \to v$ weak* in $J(u)$,
where  $v$ is a support tripotent for $J$.  Clearly 
$v u x = x$ for all $x \in J$, and hence by weak* density we have
$v u v = v$.  Also  
by the discussion above Theorem
 \ref{chopt}, $J = J(v) \cap Z$.   Now $v$ and $v^2$ are open
(since (v) implies (i) and (vi)), and Proposition \ref{3mi} (2) 
gives that $u^2 = v^2$.  
Thus $v u v = u = v$.
 \end{proof}

Parts of the last result are no doubt true for general 
open tripotents.    
We now turn to another useful way of looking at 
open tripotents.  

\begin{proposition}
\label{oppd}  Let $Z$ be a $*$-TRO, and let $A = Z + Z^2$, 
\begin{itemize} \item [(1)]  Suppose that $u$ is a selfadjoint open   tripotent in 
${\mathcal Z}(Z'')$.  Then $p = \frac{u^2 + u}{2}$ 
and $q = \frac{u^2 - u}{2}$  are open 
projections in the center of $A''$.  Moreover,
$u = p - q$ and $p q = 0$.
\item [(2)]  Suppose that 
$Z \cap Z^2 = (0)$ (this may always be ensured, by 
replacing $Z$ by a ternary $*$-isomorphic $*$-TRO).
Let $\theta : A \to A$ 
be the period 2 $*$-automorphism $\theta(z + a) = a - z$
for $a \in Z^2, z \in Z$.  Suppose that $p$ is an open
projection in the center of $A''$, such that $p q = 0$,
where $q = \theta''(p)$.
Then $u = p - q$ is a selfadjoint open tripotent in
${\mathcal Z}(Z'')$, and $q$ is also an open
projection in the center of $A''$.
\end{itemize}
   \end{proposition}

 \begin{proof}  
(1) \ It is clear that $p, q$ are projections in the center of $A''$,
and that $u = p - q$ and $p q = 0$.  
To see that $p, q$ are open, let $x_t$ be the net in 
Theorem \ref{chopt} (iii).  Then $u x_t \in Z^2$ by Lemma
\ref{uxinE} (5).
Thus $\frac{u x_t + x_t}{2}$ is a net in $A$ converging weak*
to $p$, and also $p \frac{u x_t + x_t}{2} = \frac{u x_t + x_t}{2}$.
From basic facts about Akemann's open projections it follows
that $p$ is one of these.  Similarly for $q$.
  
(2) \ Clearly $u$ is a tripotent in
${\mathcal Z}(A'')$.  Also $\theta''$ is a period 2 $*$-automorphism
of $A''$.  
Suppose that 
$(a_t)$ is an increasing net in $A$ converging weak* to $p$, with
$p a_t = a_t$.  Then $b_t = \theta(a_t)$ is a net in $A$ converging weak* to 
$\theta''(p)$, and $\theta''(p) b_t = b_t$.
Also $\theta(a_t - b_t) = b_t - a_t$.
Consequently $y_t = a_t - b_t$ is a net in $Z$
converging weak* to $u$.  So $u \in {\mathcal Z}(Z'')$.
Clearly $u^2 y_t = (p + \theta''(p)) (a_t - b_t)
= a_t - b_t = y_t$.   Moreover, 
$u y_s$ is the weak* limit of 
$$(a_t -  \theta(a_t))(a_s - \theta(a_s))
= a_t a_s + \theta(a_t a_s),$$
since $\theta(a_t) a_s = \theta(a_t) \theta''(p) p a_s = 0$.
Thus $u y_s = a_s + \theta(a_s) \geq 0$.
Thus Theorem \ref{chopt} (ii) implies that
$u = p - q$ is a selfadjoint open tripotent.   The last assertion
follows since $\theta''$ is a weak* continuous $*$-automorphism.
\end{proof}   

 The intersection of natural cones on a $*$-TRO is
again a natural cone, as may be seen by the argument
sketched above Equation (\ref{meet}).  Hence 
every family of  selfadjoint open   tripotents has 
a greatest lower bound (or `inf') amongst the selfadjoint open   tripotents.
The following fact is a little deeper, and should be important 
in future developments.  It is the `tripotent version' of 
Akemann's result that the inf of two (in this case 
central) open projections is open.

\begin{corollary} \label{meopf} Let $Z$ be a $*$-TRO, and let 
$u, v$ be two selfadjoint open   tripotents in ${\mathcal Z}(Z'')$.
Then the greatest tripotent $u \wedge v$ in ${\mathcal Z}(Z'')$
majorized by $u$ and $v$ is open, and is given by Equation
 (\ref{meet}).   Also,
 ${\mathfrak d}_u \cap {\mathfrak d}_v
= {\mathfrak d}_{u \wedge v}$, in the notation of
\ref{spl2}.
\end{corollary}  

\begin{proof}  
The last assertion is clear:
$${\mathfrak d}_u \cap {\mathfrak d}_v
= {\mathfrak c}_u \cap {\mathfrak c}_v \cap Z
= {\mathfrak c}_{u \wedge v} \cap Z =
{\mathfrak d}_{u \wedge v}.$$ 
To obtain the other fact, we appeal to Proposition \ref{oppd} (2).
We check that $u \wedge v$ (as given by Equation
 (\ref{meet})) is open.  By easy algebra,
 $$\frac{(u \wedge v)^2 + u \wedge v}{2} = 
 \frac{u^2 + u}{2} \cdot \frac{v^2 + v}{2}.$$
The latter is a product of two commuting open 
projections (by Proposition \ref{oppd} (1)), 
and hence is open.   Now the result is easy to see.       
\end{proof}   
 
{\bf Remark.}   As in JBW*-case (see e.g. \cite{ER1}), a
family of selfadjoint open central tripotents which is
bounded above by a selfadjoint open central tripotent,
has a sup.   Also,
the sum of two `orthogonal' selfadjoint open tripotents is clearly an
open tripotent.  Otherwise the `sup' of  tripotents does not make sense in 
general.

\medskip

The following is an `ordered variant' of a very useful result due to Youngson
\cite{Yng}.
We will not use this result in our paper, but it seemed worth including
 in view of the importance of Youngson's original result.

\begin{theorem} \label{orYo} Let $P$ be a completely positive
completely contractive idempotent map on a $*$-TRO $Z$.
Then $P(Z)$ is an involutive ternary system,
and $P(Z) \cap Z_+   =  P(Z_+)$ is
a natural cone on this system.   A similar assertion holds for
the matrix cones.  Thus $P(Z)$ with these matrix cones is
completely isometrically completely order isomorphic to
a $*$-TRO.
\end{theorem}

\begin{proof}  Since $P$ is positive, $P(Z_+) \subset
P(Z) \cap Z_+$.  Hence $P(Z) \cap  Z_+ = P(P(Z) \cap  Z_+) \subset P(Z_+)$.
So $P(Z_+) = P(Z) \cap Z_+$.

By Youngson's theorem \cite{Yng}, $P(Z)$ is
a ternary system with new ternary product $[Px,Py,Pz] =
P(P(x)P(y)P(z))$.  Since $[Px,Py,Pz]^* =
[Pz^*,Py^*,Px^*]$, it follows from Lemma \ref{*-TRO_char}
that $P(Z)$ is an involutive ternary system.

Since $P(Z_+) \subset Z_+$ we have
by Lemma \ref{jzp} that $P(J(Z))  \subset J(Z)$.
Hence $P$ restricts to a completely positive
completely contractive idempotent map on the C*-algebra
$J(Z)$.  By a slight variation of a well-known result of
Choi and Effros (use \cite{IAOS} Theorem 3.1 in conjunction with
\cite{CE} Lemma 3.9), $P(J(Z))$ is a
C*-algebra with respect to the new product $P(ab)$,
for $a, b \in P(J(Z))$, and the map
$x \mapsto P(x)$ from $J(Z)$ into the C*-algebra $P(J(Z))$ is
completely positive.  On the other hand
since $J(Z)$  is a ternary $*$-ideal in $Z$, if
$P(z)$ or $P(x)$
is in $P(J(Z)) \subset J(Z)$, then $P(P(x)P(y)P(z)) \subset P(J(Z))$.
We deduce that
$P(J(Z))$ is a ternary $*$-ideal in $P(Z)$.
Another part of the Choi and Effros result states that
$P(P(x)P(y)) = P(xP(y)) = P(P(x)y)$ for $x, y \in J(Z)$.
This translates to the assertion that the identity map
is a ternary $*$-morphism from
$P(J(Z))$ to $P(J(Z))$ with its C*-algebra product.
Hence $P(J(Z))$ is
a C*-ideal in $P(Z)$.    Thus by facts in the last section,
the positive cone ${\mathfrak c}$ in
$P(J(Z))$ (coming from the fact above that
$P(J(Z))$ is a
C*-algebra in a new product), is a natural cone for $P(Z)$.
Note that a representative element of ${\mathfrak c}$
is $P(a^* a)$, for $a \in P(J(Z))$, which is certainly
contained in $P(Z_+)$.   Conversely, $P(Z_+) = P(J(Z)_+) \subset
{\mathfrak c}$, by the observation above about the map $x \mapsto P(x)$
being completely positive.
Thus $P(Z_+) = {\mathfrak c}$.
We leave the `matrix cones' version as an exercise.
 \end{proof}

We end this section with a note on the title of this paper,
{\em ordered C*-modules}.  An ordered C*-module
 is an involutive C*-bimodule $Y$ over a C*-algebra $A$
(in the sense of Definition \ref{incb}),
with a given  positive cone, such that
$Y$ (with its canonical ternary product $x \langle y , z \rangle$)
is `ternary order isomorphic' to  a $*$-TRO.
Note every involutive C*-bimodule $Y$ is canonically
an involutive ternary
system (see Corollary \ref{wcor}),
and hence the notions of selfadjoint open tripotents, etc.\, make sense.
Thus by
Theorem \ref{order_char}, the ordered C*-modules
are exactly the involutive C*-bimodules $Y$,
with a positive cone of the form $\{ e \langle u | e^* \rangle : e \in Y'' \} \cap Y$,
for a 
(selfadjoint, as always) open tripotent $u$ in ${\mathcal Z}(Y'')$.
                                              
\section{Maximally ordered and unorderable $*$-TRO's}

\begin{definition}
\label{modef}   We say that a $*$-TRO $Z$ is {\em maximally ordered}  
if its given cone is maximal amongst the ordered operator space
cones on $Z$ (see \ref{odef}).  This is equivalent to saying that
every completely positive complete isometry $Z \to B$ into
a C*-algebra is a complete order injection.

We say that a $*$-TRO $Z$ is {\em unorderable}
if the only ordered operator space cone on $Z$ is the trivial one.
\end{definition}

We will need a fact about quotients of ordered ternary systems.
If $Z$ is a $*$-TRO, and if $N$ is a ternary $*$-ideal in $Z$,
then $Z/N$ is certainly an involutive ternary system (as 
remarked at the end of Section 1).  

\begin{lemma} \label{qrs}  If $J$ is a ternary $*$-ideal in 
a $*$-TRO $Z$, then
the involutive ternary system $Z/J$ possesses a
natural ordering for which the canonical
quotient ternary $*$-morphism $Z \to Z/J$ is
completely positive.
\end{lemma}
 
\begin{proof}   
We consider $Z''$ with its canonical second dual ordering.
Now $J^{\perp \perp}$ is a weak* closed 
ternary $*$-ideal in $Z''$, 
and hence equals $q Z''$ for a  central projection $q$
(see Lemma \ref{3mi}).  If $p = 1-q$ then
$(Z/J)'' \cong Z''/J^{\perp \perp}
\cong p Z''$.
We may thus identify $Z/J$ as a $*$-TRO inside the
$*$-WTRO $p Z''$.  This endows 
$Z/J$ with  natural matrix cones.  Let
$q_J : Z \to Z/J$ be the quotient ternary
$*$-morphism.  It is easy to
see that if $z \in Z_+$ then $p z \geq 0$ in $Z''$,
so that $q_J(z)$ is in the cone just defined
in $Z/J$.   A similar argument applies to matrices,
so that $q_J$ is completely positive.
\end{proof}  
            
Henceforth, whenever we refer to a natural ordering on
$Z/J$, it will be the one considered in the last Lemma.

\begin{lemma} \label{orna}  Let $Z$ be an ordered 
ternary system.   Then the given 
ordering on $Z$ is majorized by a
natural  ordering for $Z$.  In particular it follows
that a maximal ordered operator space cone on $Z$ is necessarily a
natural cone on $Z$.
\end{lemma}
 
\begin{proof}  Suppose that $T : Z \to B$ is a
completely positive complete isometry into
a C*-algebra $B$.   Let $W$ be the $*$-TRO generated by
$T(Z)$ in $B$.  By a simple result from the
sequel to this paper, there exists a canonical
surjective ternary $*$-morphism $\theta :
W \to Z$ with $\theta \circ T = I_{Z}$.  If $N = {\rm Ker} \theta$,
consider the quotient $W/N$, with the natural ordering
discussed in \ref{qrs}.  If $q_N$ is the completely positive quotient
ternary $*$-morphism $W \to W/N$,
then we obtain a surjective ternary $*$-isomorphism $\rho : Z \to
W/N$ with   $\rho = q_N \circ T$.  Note that $\rho$ is
completely positive.  Also $\rho$ induces a natural
ordering on $Z$, namely the one pulled back from
$W/N$.   The result is now clear.
\end{proof}

\begin{theorem}\label{okg}
Suppose $Z$ is an ordered ternary
system.  Then $Z$ has a maximal ordered operator space cone containing
the given one, 
and this cone is natural. 
\end{theorem}

\begin{proof}
By the previous lemma we may assume that the given order
on $Z$ is a natural order.  By Theorem \ref{order_char}, 
this order corresponds to a selfadjoint open central tripotent $u_0$.
We consider the set ${\mathcal S}$ of all selfadjoint open central tripotents 
in $E = Z''$ majorizing $u_0$, with the usual 
ordering of tripotents (see 
Proposition \ref{bij_order_corr} (1)). 
 By Theorem \ref{order_char} and \ref{orna},
we will be done if we can show that ${\mathcal S}$ has
a maximal element.  We use Zorn's lemma.  
Suppose that $(w_t)$ is an increasing net in  ${\mathcal S}$,
and let $u$ be a weak* limit point of the net.
We may assume that $w_t \to u$ weak* (by replacing the 
net by a subnet).
If $r \geq s \geq t$ then $w_r w_s w_t = w_r w_s^3 w_t
= w_s^2 w_t = w_t$.  Taking the limit over $r$ gives 
$u w_s w_t = w_t$.   Then taking the limit over $s$ 
gives
$u u w_t = w_t$.  Hence $u$ is a selfadjoint tripotent.
It is also easy to see that $u$ majorizes the net, and hence
also $u_0$.
Finally, to see that $u$ is open, we use Theorem \ref{chopt} (ii).  
Since each $w_t$ is in ${\mathfrak c}_{w_t}' \subset {\mathfrak c}_u'$,
we have $u \in {\mathfrak c}_u'$.
\end{proof}

\begin{proposition} \label{maxdo}  Let $E$ be a $*$-WTRO. 
 \begin{itemize}
\item [(1)]  If $J$ is a C*-ideal in $E$, then
$\bar{J}^{w*}$ is a weak* closed ternary $*$-ideal in $E$ which
is ternary $*$-isomorphic to a W*-algebra.
\item [(2)]  Every ordered operator space cone on 
$E$ is contained in a natural dual cone.
\item [(3)]  The maximal ordered operator space cones 
on $E$ coincide with the 
maximal natural dual cones on $E$.
\end{itemize}  \end{proposition}

\begin{proof}  (1) \ It is easy to see that $\bar{J}^{w*}$ is a weak* 
closed ternary $*$-ideal in $E$.
Suppose that $A$ is a C*-algebra,
and that $\varphi : A \to J$ is a ternary $*$-isomorphism.
Let $(e_t)$ be an increasing approximate identity for 
$A$. We may assume that $\varphi(e_t) \to u \in \bar{J}^{w*}$
in the weak* topology.   Clearly $u$ is selfadjoint,
and $u \varphi(a) \varphi(b) = \lim \varphi(e_t ab) 
= \varphi(ab)$ for all $a, b \in A$.
Taking $a  = e_t$, and taking the limit, yields 
$u^2 y = y$ for all $y \in J$.   Hence $u^2 y = y$ for all $y \in 
\bar{J}^{w*}$.  Thus $u$ is a tripotent in $\bar{J}^{w*}$.
Clearly the above argument also gives
$\varphi(a) \varphi(b) u = \varphi(ab)$ for all $a, b \in A$,
so that $u$ is in the center of $\bar{J}^{w*}$.  It follows that
$u \bar{J}^{w*}$ is a W*-subalgebra of the W*-algebra 
$\overline{E^2}^{w*}$.  Thus by Lemma 
\ref{w*-ideal_char} we deduce that
$u \in {\mathcal Z}(E)$, and that $\bar{J}^{w*} = J(u)$.

(2) \ Any ordered operator space cone 
on $E$ is contained in a natural cone,
by Lemma \ref{orna}.  
Suppose that ${\mathfrak C}$ is a natural cone
on $E$; thus $(E,{\mathfrak C})$ 
is ternary order isomorphic to a $*$-TRO $Z$.
Then $E$ has a ternary $*$-ideal $J$ which is ternary $*$-isomorphic to 
$J(Z)$.  Thus by (1), $E$ contains a 
weak* closed ternary $*$-ideal $J'$ containing $J$, which
is ternary $*$-isomorphic to a W*-algebra.  By Lemma
\ref{w*-ideal_char}, $J' = J(w)$ 
for a selfadjoint tripotent $w$ in 
${\mathcal Z}(E)$.   Also ${\mathfrak C} \subset
{\mathfrak c}_w$, since $J$ is a
C*-subalgebra of $J(w)$.  The result is completed by
an appeal to  Proposition \ref{bij_order_corr} (4).
Similarly for the matrix cones.

(3) \ Follows from (2).   \end{proof}

The following result follows immediately from the above and Corollary
 \ref{maxcnd}.

\begin{corollary}\label{chmo}  A $*$-WTRO $E$  
  is  maximally ordered if and only if $E$ is completely order isomorphic via
a ternary $*$-isomorphism to $M \oplus^\infty F$, where $M$ is a W*-algebra
and $F$ is a $*$-WTRO which is unorderable.  Also,
 $E$ is unorderable if and only if 
${\mathcal Z}(E) = \{0\}$.   \end{corollary}                          

\begin{corollary} \label{strot}
Let $Z$ be an involutive ternary system. Then
$Z$ is unorderable if and only if $Z$ contains no nonzero
C*-ideals, and if and only if there are no 
nonzero selfadjoint open tripotents in ${\mathcal Z}(Z'')$.
\end{corollary}

\begin{corollary} \label{ze0}  If $Z$ is a $*$-TRO,
and if  ${\mathcal Z}(Z'') = \{ 0 \}$, then $Z$ is
unorderable; and also ${\mathcal Z}(Z) = \{ 0 \}$.
\end{corollary}
 
We do not have a characterization of 
maximally ordered $*$-TRO's which is quite as tidy as
in the $*$-WTRO case.  In fact the situation here
seems quite complicated, as the examples towards the end
of Section \ref{ilb} show.
Here are a couple of partial results:  
 
\begin{corollary}\label{sdmao}
Let  $Z$ be a $*$-TRO such that $Z''$ with its canonical
ordering is maximally ordered
(or equivalently such that
the selfadjoint open central tripotent corresponding to the given
ordering is unitary).  Then:
\begin{itemize} \item [(1)]  $Z$ is maximally ordered. 
\item [(2)]  The involutive ternary system
$Z/J(Z)$ is unorderable.
\end{itemize}
In particular a $*$-TRO $Z$
is maximally ordered if $Z''$ is a C*-algebra.
\end{corollary}      

\begin{proof}
(1) \ If $T$ is a completely positive complete isometry $Z \to B$ into
a C*-algebra then
$T''$ is completely positive by Proposition \ref{jinjs}.
By hypothesis, $T''$ is a complete order
embedding, and therefore so is $T$, by Proposition \ref{jinjs} again.

(2) \ 
If $Z/J(Z)$ had a
nontrivial natural ordering, then so does
$(Z/J(Z))'' \cong Z''/J(Z)^{\perp \perp}$.
However since $J(Z)^{\perp \perp} = J(u)$ for
a selfadjoint unitary $u$ in ${\mathcal Z}(Z'')$,
it follows that $Z''/J(Z)^{\perp \perp}
\cong \{ e - u e u : e \in Z'' \}$.  Hence the latter
involutive ternary system has a
nontrivial natural dual ordering, and hence has
a nonzero selfadjoint tripotent $w
\in {\mathcal Z}(Z'')$ (by Theorem \ref{w*-natcone_char}).
Since $w + u \geq u$, this contradicts
Proposition \ref{bij_order_corr} (2).

Finally note that C*-algebras are maximally ordered, as one can see
for example by combining (1) with Lemma \ref{orna}
(or see \cite{GILO,Wi}).
\end{proof}
 
{\bf Remarks.}  1)  \   The converse of Corollary \ref{sdmao} (1) is false, as we show
in  an example towards the end of Section
\ref{ilb}. That is, $Z$ may be maximally ordered without $Z''$
(with its canonical cone) being maximally ordered.

2) \    It is not true in general 
that for a maximally  ordered $*$-TRO $Z$,
$Z/J(Z)$ is necessarily unorderable.  Counterexamples may be found
amongst the examples considered towards the end of 
Section \ref{ilb}, together with the fact
(proved above Proposition \ref{trsidlb2} )
that the quotient of a commutative involutive
C*-bimodule (in the sense of \ref{incb})
by a ternary $*$-ideal is again a
commutative involutive
C*-bimodule.   However we do have the following result (which also
gives another proof of Corollary \ref{sdmao} (1)):

\begin{proposition} \label{zoj}  If $Z$ is a $*$-TRO
such that $Z/J(Z)$ is unorderable, then $Z$ is maximally ordered.
\end{proposition}

\begin{proof}
Let $v$ be the open tripotent corresponding to the given ordering
on $Z$.  Let $E = Z''$ as usual,
then $J(E) = J(v) = J(Z)^{\perp \perp}$ (see 
Proposition \ref{jinjs1}).  If $Z$ is not 
maximally ordered, then there exists an open tripotent $u \geq v$,
$u \neq v$.  We identify the $*$-TRO 
$\{ x - v x v : x \in E \}$ with $(Z/J(Z))''$ via the 
canonical isomorphisms
$$\{ x - v x v : x \in E \} \cong E/J(v) \cong (Z/J(Z))''.$$
Then $u - v$ is a nonzero central tripotent in this $*$-TRO; if we can 
show that it is also open we will have the desired contradiction,
by Corollary \ref{strot}.    If $(x_t)$ is as in Theorem 
\ref{chopt} (iii), then it is easy to check using
Theorem
\ref{chopt} (iv) that $x_t + J(Z) = x_t + J_v(Z)$ converges
in the weak* topology of $(Z/J(Z))''$ to $u - v$.  
We conclude with an appeal to Theorem
\ref{chopt} (ii).  \end{proof}

We now turn to a variant of the topological boundary of an
open set.  Let $A$ be a C*-algebra, and let $p$ and $q$ be respectively
an open and a closed central projection in $A''$.   We say that
$q$ {\em is contained in the boundary of} $p$, if
$p q = 0$, and if $r p \neq 0$ whenever $r$ is an open
central projection in $A''$ such that $r q \neq 0$.

\begin{proposition} \label{chbdy}  Let $Z$ be a $*$-TRO, let
$u$ be a selfadjoint open tripotent in ${\mathcal Z}(Z'')$, and define
$p$ and $q$ as in Proposition \ref{oppd} (1).  Write $1$ for
the identity of ${\mathcal Z}(Z'')^2$ (see e.g.\ the
proof of Lemma \ref{z(e)w*} (1)).
Suppose that $1 - u^2$ is contained in the boundary of
both $p$ and $q$.   Then
$u$ is a maximal selfadjoint open tripotent in ${\mathcal Z}(Z'')$.
Hence the corresponding cone
${\mathfrak d}_u$ is a maximal cone for $Z$.
\end{proposition}

\begin{proof}     Suppose that $v \geq u$, but $v \neq u$.
Then $u^2 \neq v^2$.  Hence either $\frac{v^2 + v}{2} \neq \frac{u^2 + u}{2}$
or $\frac{v^2 - v}{2} \neq \frac{u^2 - u}{2}$.
Assume the former (the other case is similar).
Thus $\frac{v^2 + v}{2} (1-u^2) \neq 0$.  By hypothesis,
$\frac{v^2 + v}{2} \frac{u^2 - u}{2} \neq 0$.
However $(v^2 + v)(u^2 - u) = u^2 + u - u - u^2 = 0$,
a contradiction.
 \end{proof}

{\bf Remark.}  We imagine that a modification of the ideas in
the last Corollary yields a {\em characterization} of 
maximal selfadjoint open tripotents, and therefore also of maximal cones.
Indeed in the commutative case the `boundary' hypothesis
in Proposition \ref{chbdy} is necessary and sufficient (see
the later Corollary \ref{ntoc}). 
In the general case this converse requires further investigation.

\medskip
 
{\bf Example.}  
If $M$ is a finite dimensional W*-algebra, it is easy to see
that the cones on $M$ may be characterized as follows.  If 
the center of $M$ is $n$ dimensional,
the span of $n$ minimal central projections $\{ p_1 ,
\cdots , p_n \}$, then there are $3^n$ possible
natural cones on $M$, namely the product of the natural cone of $M$, with each
of the $3^n$ selfadjoint tripotents $(\alpha_1 p_1) \oplus
(\alpha_2 p_2)  \cdots (\alpha_n p_n)$, where $\alpha_k \in \{ -1 , 0 , 1 \}$.
There are $2^n$  maximal ordered operator space cones on $X$, namely the cones
in (2) above, but with $\alpha_k \in \{ -1 , 1 \}$.

\medskip
 
{\bf Example.}   If $E$ is a finite dimensional involutive ternary system,
then as above we may write $E \cong M \oplus^\infty F$, where
$M$ is a finite dimensional W*-algebra, and $F$ is 
an unorderable finite dimensional involutive ternary system.
The possible natural cones on $E$, are exactly those in the last example.
 
\medskip
      
In fact it is not hard to characterize the unorderable finite dimensional 
involutive ternary systems.  They are precisely the 
$*$-TRO's $e N (1-e) + (1-e) N e$, for
a finite dimensional W*-algebra $N$ and a projection $e \in N$. 
There is a similar characterization valid for all `type I' unorderable 
involutive dual ternary systems.
However in general it is not true that 
general (non type I) unorderable involutive dual ternary systems
are of the form $e N (1-e) + (1-e) N e$, for
a W*-algebra $N$ and a projection $e \in N$.
We will discuss these matters elsewhere.

We will also discuss unitizations of $*$-TRO's elsewhere.

\section{Examples: Involutive line bundles} \label{ilb}
 
In this section we consider a special class of
involutive ternary systems, namely the $*$-TRO's $Z$ in  
$C(K)$ spaces.  We begin by remarking that in this case the `center'
${\mathcal Z}(Z)$ (defined in Section \ref{ICB})
is simply $Z$, and the `center' of
$Z''$ is simply $Z''$.  It follows from the earlier
theory that $Z''$ is ternary $*$-isomorphic to an
$L^\infty$ space, a fact that we will use
(sometimes silently) below.  Also, one might guess
by analogy to Corollary \ref{chmo}, that
such involutive ternary systems
always have nontrivial orderings, and in fact we shall see
that this is the case.

The $*$-TRO's $Z$ in
$C(K)$ spaces may be given several abstract 
characterizations.
 For example they are the 
commutative involutive C*-bimodules of Definition
\ref{incb}.
They may also be viewed as being
the space of sections vanishing at infinity for certain
line bundles.
We omit the proof of these assertions, and instead
give the following abstract
characterization related to the notion of `$C_\sigma$ spaces'
(see e.g.\ \cite{FRc,Kaup3}).
Namely, suppose that $\Omega$ is a locally compact Hausdorff 
space, and that $\tau : \Omega \to \Omega$ is a homeomorphism
with $\tau \circ \tau = I_\Omega$.   This corresponds to
a period 2 automorphism of $C_0(\Omega)$.  Let  $W$ be 
the corresponding $*$-TRO, namely
\begin{equation} \label{cilb}
W = \{ f \in C_0(\Omega) : f \circ \tau = - f \}.
\end{equation}
Then $W$ is a $*$-TRO in 
$C_0(\Omega)$.  In fact every $*$-TRO in
a commutative
C*-algebra arises in this way.  This is essentially well
known (see e.g.\ \cite{SOJ,Kaup3} and references therein), 
but for completeness we will include a proof which includes
details we will need later.   

\begin{theorem} \label{evcl}  If $Z$ is a $*$-TRO in a commutative
C*-algebra, then $Z$ is ternary $*$-isomorphic to 
a ternary system $W$ of the form in (\ref{cilb}),
where $\Omega$ is a locally compact subspace of 
$Ball(Z')$ with the weak* topology, and $\tau$ is
simply `change of sign': namely $\tau(\varphi) = - \varphi$.
Indeed $\Omega$ may be taken to be
$ext(Ball(Z'_{sa}))$.   \end{theorem}

\begin{proof}  
Let $Z$ be a closed ternary $*$-subsystem  of a commutative 
C*-algebra.   
We recall two facts that may be found in
\cite{FRc}: firstly, 
$ext(Ball(Z'))$ is
a locally compact space with the weak* topology,
and secondly that any $\psi \in ext(Ball(Z'))$ is
a ternary morphism on $Z$.  Now 
$Z'_{sa}$, the set of selfadjoint functionals in $Z'$,
 is closed in $Z'$ in the weak* topology,
so that $Ball(Z'_{sa})$ has plenty of nonzero extreme points
by the Krein-Milman theorem.
Set $\Omega = ext(Ball(Z'_{sa}))$, we claim that this
is a subset of $ext(Ball(Z'))$.  To see this
let $\psi \in ext(Ball(Z'_{sa}))$, 
and let $\varphi_1,
\varphi_2 \in Ball(Z')$ with
$\psi = \frac{\varphi_1 + \varphi_2}{2}$.  We may write
$\varphi_k = \rho_k + i \sigma_k$, where 
$\rho_k, \sigma_k \in Ball(Z'_{sa})$.  
Applying $\psi$ to elements of
$Z_{sa}$ and taking real parts 
shows that $\sigma_1 + \sigma_2 = 0$ on $Z_{sa}$, and hence
$\sigma_1 + \sigma_2 = 0$ on $Z$.  Thus
$\psi = \frac{\rho_1 + \rho_2}{2}$, so that $\rho_1 = \rho_2 = \psi$.
We will use the fact that $Z''$ is ternary $*$-isomorphic
to a commutative W*-algebra
$L^\infty(\Omega,\mu)$.    In particular note that then
$Z'$ is isometrically $*$-isomorphic to 
$L^1(\Omega,\mu)$.  However if $f, g \in 
Ball(L^1(\Omega,\Rdb))$ with 
$\Vert f + i g \Vert_1 \leq 1 = \Vert f \Vert_1$, then
it is easy to see that $g = 0$.   Thus it follows 
that $\varphi_1 = \psi + i \sigma_1 = \psi$, and 
similarly $\varphi_2 = \psi$.    Thus we have proved the claim.
Since $Z'_{sa}$ is closed in $Z'$ in the weak* topology,
it is easy to check
that $\Omega$ is also a closed subset of $ext(Ball(Z'))$,
and therefore is locally compact in the weak* topology.  

Let 
\begin{equation} \label{defW} 
W = \{ f \in C_0(\Omega) : f(-\psi) = - f(\psi)
\; \text{for all} \; \psi \in \Omega \}.
\end{equation} 
 Then $W$ is a space of the form (\ref{cilb}).
Define $\Phi : Z \to W$ by $\Phi(z)(\psi) = \psi(z)$,
for $z \in Z, \psi \in \Omega$.
  This is a ternary $*$-morphism by remarks in the last paragraph, 
and is easily checked to be
 1-1.  Hence $\Phi$ is an isometry.   We claim 
that $\Phi$ is surjective; this follows from a `Stone-Weierstrass
theorem for line bundles' such as Theorem 4.20 in \cite{BShi}.   
More specifically, note that
$B = \{ f \in C_0(\Omega) : 
f(-\psi) = f(\psi) \; \text{for all} \; \psi \in \Omega 
\}$ is a commutative C*-algebra.   One may define an equivalence
relation $\equiv$ on $\Omega$, by identifying $-\psi$ and $\psi$.
Then  $S = \Omega/\equiv$ is a locally compact Hausdorff 
space, and $B \cong C_0(S)$ as $\cas$.  It is easy to
see that $W^2$ strongly separates points of $S$, so that
by the Stone-Weierstrass theorem $W^2 = B \cong C_0(S)$.
Also $\Phi(Z)^2$ strongly separates points of $S$, so that
by Theorem 4.20 in \cite{BShi} it follows that 
$\Phi$ is surjective.
\end{proof}

The reader may notice that some of the arguments/results below may be 
replaced by references to facts in 
in the earlier sections.  We have chosen not to do this when a simple
self-contained argument came to hand.  

A subset $C$ of $\Omega$ will be called
{\em symmetric} if
$C = -C$, and {\em antisymmetric} if $C \cap (-C) = 
\emptyset$.   
Given a  closed symmetric subset $C$ of $\Omega$,
and if $W$ is as in the last proof, then 
we will write 
\begin{equation} \label{defFC}
F_C = \{  f \in W : f(x)  = 0 \; \text{for all}
\;  x \in C \}.  \end{equation}    

If $S$ is as in the last proof, then we saw that
$W^2 \cong C_0(S)$ $*$-isomorphically.  Also, the 
open projections in $C_0(S)''$ correspond to the
open sets in $S$, and hence to the symmetric
open sets in $\Omega$.  This together with Lemma
\ref{3mi} yields:


\begin{proposition} \label{trsidlb}  A subspace $N$ 
of $W$ is a ternary $*$-ideal of $W$ if and only 
if there is a closed symmetric subset $C$ of $\Omega$ 
such that $N = F_C$.
Moreover such a set $C$ is uniquely determined by $N$.     
\end{proposition}     

{\bf Remark.}  If $W$ is as defined in (\ref{defW}), and if $C$ is a
closed symmetric subset of $\Omega$, define a
ternary $*$-morphism $W \to \{ f \in C_b(C) :
f(-\psi) = - f(\psi) \; \text{for all} \;  
\psi \in C \}$ by $f \mapsto f_{|_C}$.   The 
kernel of this ternary morphism is $F_C$, so that
$W/F_C$ is ternary $*$-isomorphic to a
$*$-TRO in $C_b(C)$.  This shows that
the class we are investigating in this section, 
namely the $*$-TRO's in  commutative
$\cas$, is closed under quotients by ternary $*$-ideals.

\smallskip 
 
It is easy to see that C*-ideals in $W$ are necessarily
ternary $*$-isomorphic to commutative C*-algebras.  We will use this
fact several times below.
 
\begin{proposition} \label{trsidlb2}   Let 
$\Omega$ be as defined in Theorem \ref{evcl},
and let  $U \subset \Omega$ be an open antisymmetric
subset, and let $C = 
\Omega \setminus (U \cup (-U))$.  Then
$F_C$ is a C*-ideal.
Conversely,
if
$C$ is a closed subset $C$ of $\Omega$ with
$C = -C$,
such that
$F_C$ is a
C*-ideal,
then there is an open set $U \subset \Omega$ with
$U \cap (-U) = \emptyset$, and $U \cup (-U) = C^c$.
\end{proposition}
 
\begin{proof}   
Given $U$ with the
asserted properties, define 
$\theta : C_0(U) \to F_C$ by $\theta(f)(x) = f(x)$ 
if $x \in U$, $\theta(f)(x) = -f(-x)$
if $x \in (-U)$, and $\theta(f)(x) = 0$ otherwise.
Note that $\theta(f)$ is continuous on $U$, and therefore
also on $-U$.  In the interior of $C$ the function
$\theta(f)$ is clearly
continuous.  Finally if $x$ is in the boundary of $C$,
and if $\epsilon > 0$ is given then there is a compact
$K \subset U$ with $|f| < \epsilon$ on $U \setminus K$.
Choose an open symmetric set $V$ containing $x$ which
does not intersect $K \cup (-K)$. On $V$ it is 
clear that $|\theta(f)| < \epsilon$.  Thus $\theta(f)$ is continuous 
at $x$.  Indeed it is clear that  $\theta(f) \in C_0(\Omega)$,
so that  $\theta(f) \in F_C$.   Clearly $\theta$ is
a 1-1 ternary $*$-morphism.   Given $g \in F_C$ let 
$f $ be the restriction of $g$ to $U$.  If $\epsilon > 0$ is given
let  $K = \{ w \in \Omega : |g(w)| \geq \epsilon \}$,
this is compact and is a symmetric subset of $U \cup (-U)$.
Thus $K \cap U$ is a compact subset of $U$, so that
$f \in C_0(U)$.  Clearly $\theta(f) = g$. 

Conversely, suppose that $B$ is a commutative C*-algebra,
and that 
$\theta : B \to W \subset C_0(\Omega)$ is a 1-1 ternary $*$-morphism onto
$N = F_C$.  Let $(e_t)$ be an increasing
approximate identity for $B$.  Then the weak* limit $v$ of $(\theta(e_t))$
is the support tripotent for $N$ (see the discussion before
Theorem \ref{chopt}).  For $\omega \in \Omega$ let 
$\delta_\omega$ be `evaluation at $\omega$'. This
yields  a canonical map
$\Omega \to W : \omega \mapsto \delta_\omega$.   
Define  $h(\omega) = v(\delta_\omega)$ for $\omega \in \Omega$.
Clearly $h$ is a real valued function on $\Omega$, and  
$h(\psi) = \lim_t \theta(e_t)(\psi)$ for all
$\psi \in \Omega$.
If $g = \theta(b)$ for $b \in B$, then we have 
\begin{equation} \label{eqnd}
h(\psi) g(\psi)^2 = \lim_t \theta(e_t)(\psi) g(\psi)^2 =
\lim_t \theta(e_t b^2)(\psi) = \theta(b^2)(\psi).
\end{equation} 
 Similarly, $h^2 g = g$
for all $g \in N$.
Thus $h(\psi) = 1$ or 
$-1$ for every $\psi \notin C$.   Let $U = h^{-1}({1})$,
then $-U = h^{-1}({-1})$, and these are disjoint subsets of
$\Omega \setminus C$.   If $\psi \notin C$ and 
$g \in N$ with $|g(\psi)| = \alpha \neq 0$, let
$V = |g|^{-1}((\alpha/2,\infty))$ a symmetric open set in
$\Omega \setminus C$, containing $\psi$.  Now by (\ref{eqnd}), $h g^2 = 
\theta(b^2)$ is continuous on $V$,
so that $h$ is continuous on $V$. 
Hence $h$ is continuous on $\Omega \setminus C$,
and so $U$ and $-U$ are open.
 \end{proof}

\begin{definition} Given an open antisymmetric set $U \subset \Omega$ as 
in the last proposition, define ${\mathfrak c}_U$ 
to be the cone in
$W$ consisting
of those functions $f \in W$ with $f(x) \geq 0$ if $x \in U$,
and $f(x) = 0$ for $x \in C = \Omega \setminus (U \cup (-U))$.
(Necessarily then $f \leq 0$ on $-U$.)   
\end{definition}

{\bf Remark.}  In the last proof we saw that the support
tripotent $v$ of a C*-ideal `restricts' to a Borel function
$h$ on $\Omega$.    There is a canonical map
$\rho : Bo(\Omega) \to C_0(\Omega)''$, where
$Bo(\Omega)$ is the C*-algebra of bounded Borel measurable
functions on $\Omega$.  We claim that $\rho(h) = v$.  To see
this, it suffices (by Krein-Milman) that
$\rho(h)(\psi) = v(\psi)$ for all $\psi \in ext(Ball(W'_{sa}))$.
However any such $\psi$ is just one of the $\delta_\omega$ in
the last proof, and then the desired relation follows from
the definition of $h$.    As a corollary one may deduce that
${\mathfrak c}_U = {\mathfrak c}_v \cap W$.
This fact could be used to shorten some of the arguments below.

\begin{lemma}  \label{ceq}  If $U_1, U_2$ are two
open antisymmetric subsets of $\Omega$,
then ${\mathfrak c}_{U_1} = {\mathfrak c}_{U_2}$ if and only
if $U_1 = U_2$.  \end{lemma}  

\begin{proof}   The one direction is trivial.
If $U_1 \neq U_2$ and if $x \in U_1 \setminus U_2$
then by Urysohn's lemma we can choose $f \in C_0(U_1)$
with $f \geq 0$ and $f(x) = 1$.  By the first part of the proof of 
\ref{trsidlb2} we obtain a function $g \in {\mathfrak c}_{U_1}$
with $g(x) = 1$.  Since $x \notin U_2$ we have $g(x) \leq 0$.
Thus ${\mathfrak c}_{U_1} \neq {\mathfrak c}_{U_2}$.    
 \end{proof}    

\begin{corollary} \label{putic}  
\begin{itemize} \item [(1)]  If $U$ 
is an  open antisymmetric set $U \subset \Omega$,
then ${\mathfrak c}_U$ is a
natural cone.
\item [(2)]  Given  
a natural cone ${\mathfrak c}$ in $W$,
then there is a unique open antisymmetric set $U \subset \Omega$ 
such that ${\mathfrak c} = {\mathfrak c}_U$.
Indeed $U = \cup_{g \in {\mathfrak c}} \;
g^{-1}((0,\infty))$. 
\end{itemize} \end{corollary} 

\begin{proof}  Part (1) is obvious from the first part of the proof of
Proposition \ref{trsidlb2}.
To see (2), first note that if ${\mathfrak c}$ is
a natural cone in $W$, then $N = \; $Span$({\mathfrak c})$
is a C*-ideal.  By Propositions
\ref{trsidlb} and \ref{trsidlb2}, we obtain the
associated open set $U$.
Let $\theta$ be as in the (second part of the) last proof.
If $\pi : B \to N^2$ is the canonical $*$-isomorphism
associated with $\theta$ (see Proposition \ref{prtm} (6)), 
then $\theta(b) = v \pi(b)$ for all $b \in B$.  
Indeed
$$\theta(b)(\psi) = \lim_t \theta(e_t^2 b)(\psi)
= \lim_t \theta(e_t)(\psi) \lim_t (\theta(e_t) \theta(b))(\psi)
= v(\psi) \pi(b)(\psi).$$
Thus since $b \geq 0$ if and only if 
$\pi(b) \geq 0$, and since 
$v = 1$ on the open set $U$,
it follows that a function $g$ in $F_C$
is in ${\mathfrak c}$ 
(i.e.\  of the form $\theta(b)$  for a $b \in B_+$) if and only if 
$g \in {\mathfrak c}_U$.
Also note that
$x \in U$ if and only if $g(x) > 0$ for some 
$g \in \theta(B_+)$.  To see this note that
if $x \in U$ then $g(x) \geq 0$ for all
$g \in \theta(B_+)$ by the above.  If 
$g(x) = 0$ for all
$g \in \theta(B_+)$, then $g(x) = 0$ for all
$g \in \theta(B) = N$, which is impossible.
Conversely. if $g(x) > 0$ for some
$g \in \theta(B_+)$ then certainly $x \notin C$ since
$g \in F_C$.   If $x \in (-U)$ then $-x \in U$, so that
$g(-x) = - g(x) \geq 0$, which is impossible.  So
$x \in U$. \end{proof}

\begin{lemma} \label{odcon}
Suppose that we have two natural cones ${\mathfrak c}_1$ 
and ${\mathfrak c}_2$ on 
$Z$, or equivalently on $W$, 
and that $U_1$ and $U_2$ are the two corresponding open antisymmetric
subsets of $\Omega$. 
Then
$U_1 \subset U_2$ if and only if 
${\mathfrak c}_1 \subset {\mathfrak c}_2$.
\end{lemma} 

\begin{proof}   The key point is that ${\mathfrak c}_k
= {\mathfrak c}_{U_k}$.   
Thus one direction is obvious: if $U_1 \subset U_2$.  
and if  
$g \geq 0$ on $U_1$, and $g = 0$ on $\Omega \setminus
(U_1 \cup (-U_1))$, then $g \geq 0$ on $U_2$.
The converse direction follows from the relation
$U_k = \cup_{g \in {\mathfrak c}_k} \;
g^{-1}((0,\infty))$ from Corollary \ref{putic}. 
\end{proof}
  
\begin{corollary}   \label{nto}   Let
$Z$ be a nontrivial (i.e.\  nonzero)
commutative involutive ternary system.
Then $Z$ possesses a nontrivial  natural ordering.
\end{corollary}

\begin{proof}   We may assume by
\ref{evcl} that $Z = W$, where 
$W$ is as defined in (\ref{defW}).  Let $\psi \in \Omega$
and pick $x \in Z$ such that $\psi(x) = 1$.
Then  $U = \{ \varphi \in \Omega : 
|\varphi(x) - 1| < \frac{1}{2} \}$ is an 
open subset of $\Omega$ which does not intersect
$-U$.
By Proposition \ref{trsidlb2} we obtain a 
nontrivial natural ordering on $Z$.    
\end{proof}

Putting together Corollary \ref{nto} and Theorem
\ref{okg}  we have:

\begin{corollary}   \label{lbmaco}  Let
$Z$ be a nontrivial (i.e.\  nonzero)
commutative involutive ternary system.
Then $Z$ possesses a nontrivial
maximal ordering, and
this  ordering is a natural one.
\end{corollary}

\begin{corollary}   \label{ntoc}   For a
nontrivial (i.e.\  nonzero)
commutative involutive ternary system $Z$,
let $W$ be as defined in (\ref{defW})
 (so that $Z \cong W \subset 
C_0(\Omega)$, where $\Omega \subset Ball(Z'_{sa})$).
If $U$ is an open subset of $\Omega$ with
$U \cap (-U) = \emptyset$, then the following are equivalent:
\begin{itemize}
\item [(i)]  ${\mathfrak c}_U$ (the functions in $W$ which 
are nonnegative on $U$) is a  maximal natural cone of $W$,
\item [(ii)] $\Omega \setminus (U \cup (-U)) = Bdy(U)$,
\item [(iii)]  $\Omega \setminus (U \cup (-U))  = Bdy(-U)$,
\item [(iv)]   $Bdy(U) = Bdy(-U)$ and 
$\Omega \setminus (U \cup (-U))$ has no interior points,
\item [(v)]  there is no larger open subset $U'$ of
$\Omega$ containing $U$ with $U' \cap (-U')= \emptyset$.
 \end{itemize}    
Any maximal cone of $W$ is of the form in {\rm (i)}.
  \end{corollary} 
   
\begin{proof}  Let $C = \Omega \setminus (U \cup (-U))$.
Both $Bdy(U)$ and $Bdy(-U)$ are contained in $C$.

(v) $\Rightarrow$ (iii) \  Suppose that  $C \neq Bdy(-U)$.
Then $C$ contains a point $\varphi$ which is an
interior point for $C \cup U$.   Suppose that 
$\varphi(z) = 1$ for some $z \in Z$.
Then if $\varphi
 \in V \subset C \cup U$ with $V$ open, let
$O = 
V \cap \{ \psi \in \Omega :  |\psi(z) - 1| < \frac{1}{2} \}$.    This is an
open set in $\Omega$ containing $\varphi$, and
$(O \cup U) \cap (-O \cup -U) = \emptyset$.         
To see this note that it is quite obvious
that $O \cap (-O) = \emptyset$.  If $\psi \in O \cap
(-U)$ then $\psi \in C$.   But 
$C \cap
(-U)  = \emptyset$.  Thus $O \cap
(-U) = \emptyset$, and also $(-O) \cap U = \emptyset$.  

(iii)  $\Rightarrow$ (v) \  Suppose that  $C = Bdy(-U)$.
If $U'$ is a larger open set as in (v), and if 
$x \in U' \setminus U$, then $x \in C = Bdy(-U)$.
Thus $U' \cap
(-U) \neq \emptyset$, which is impossible. 
So $U = U'$.  

 (ii) $\Leftrightarrow$ (iii)
 $\Leftrightarrow$ (iv) \ Easy, since $C = -C$.

(v) $\Rightarrow$ (i) \  Suppose that
we had a natural cone containing ${\mathfrak c}$.
Then these two
orderings must be the same, by Lemma \ref{odcon} and the 
remark before it, and
the hypothesis (v). 

(i) $\Rightarrow$ (v) \  $U \subset U'$ implies that
${\mathfrak c}_U \subset {\mathfrak c}_{U'}$.
By hypothesis ${\mathfrak c}_U = {\mathfrak c}_{U'}$,
so that $U = U'$ by Lemma \ref{ceq}.    

The last statement uses also 
Lemma \ref{orna} and 
Corollary \ref{putic}.  \end{proof} 

From Corollary \ref{ntoc} and Lemma \ref{ceq},
we see that there is 
a bijection between maximal cones (which we know from an earlier
result are always natural) on a
commutative involutive ternary system, and 
the class of open antisymmetric subsets $U$ 
characterized in Corollary
\ref{ntoc}.  Indeed we may call such a subset $U$ an 
`antipodal coloring' of $\Omega$,
and then maximal cones are in
1-1 correspondence with such `antipodal colorings'.
From this correspondence
it is easy to construct interesting very explicit examples
of maximal cones on  commutative involutive ternary systems.
The example we shall consider in the remainder of
this section
is as follows.
Let $S^2$ be the unit sphere, and $Z$ the $*$-TRO
$\{ f \in C(S^2) : f(-x) = - f(x) \}$.   This is clearly
a trivially ordered $*$-TRO in $C(S^2)$.
By an `antipodal coloring of the sphere' we mean an open
subset $U$  of the sphere (called blue), which does not intersect
$-U$ (called red), such that the boundary of $U$ is the boundary of $(-U)$,
and this latter set has no interior.
Thinking about such colorings of the sphere it is clear that
there is a rich profusion of them that are quite
different topologically. One may also, if one wishes,
choose the coloring
so that 
the measure of $C = \Omega \setminus (U \cap (-U))$ is positive.

\vspace{ 3 mm}

{\bf Example.}   
In this example we show that given a
natural cone ${\mathfrak c}$ on a $*$-TRO $Z$,
there may exist
many distinct maximal orderings on the second dual $E =
Z''$, whose restriction
to $Z$ is ${\mathfrak c}$.
This can be done even when ${\mathfrak c}$ is $(0)$ or is
a maximal ordering on $Z$.       
We will use the fact
from Proposition \ref{bij_order_corr}
that any selfadjoint unitary $u$ in
${\mathcal Z}(E)$ gives a maximal cone ${\mathfrak c}$  on $E$
such that $z \in {\mathfrak c}$ if and only if $u z \geq 0$.  
 Let $Z = \{ f \in C(S^2) : f(-x) = - f(x) \}$, and let
$H$ be the open upper hemisphere of $S^2$.
Let $P$ and $Q$ be disjoint Borel sets which
together 
partition $H$, each of which is dense in $H$.
Let $\Tdb$ be the `equator' in $S^2$.
Let $v$ be the function
from $H$ to $\{-1, 1\}$ which is $1$ on $P$ and $-1$ on $Q$.
Then $v$ is a unitary in the C*-algebra $Bo(H)$
of bounded Borel functions on $H$.  Let $\theta :
Bo(H) \to C_0(H)''$ be the canonical unital $*$-homomorphism.
There is a canonical
ternary $*$-morphism $\nu : C_0(H) \to Z$ as in the proof
of \ref{trsidlb2}.
 Then $\psi = \nu'' : 
C_0(H)'' \to Z''$ is a 1-1 ternary $*$-morphism.
Let $\rho = \psi \circ \theta$ and $V = \rho(v)$.
Then $V$ is a selfadjoint tripotent in $Z''$.
We may choose, as in Lemma \ref{z(e)w*} (3), a 
selfadjoint unitary $u$ in $Z''$ with $u \geq V$. 
As we saw in Proposition
\ref{bij_order_corr},  we may endow $Z''$  with
a maximal cone ${\mathfrak c}_u$
(which depends on $P$).
We claim that the restriction of ${\mathfrak c}$  to
$Z$ is trivial.  Indeed if $f \in Z$ with $f u \geq 0$, 
then $f u V^2 = f V \geq 0$.
If $\omega \in H$, and if $\delta_\omega$ is `evaluation
at $\omega$' (which is a character of $C(S^2)$),
then we have $(f V)(\delta_\omega) = f(\omega) V(\delta_\omega)
\geq 0$.  However $V(\delta_\omega) = 1$ if $\omega \in P$, because
$$V(\delta_\omega) = \nu''(\theta(v))(\delta_\omega)
= \theta(v)(\nu'(\delta_\omega))) = \theta(v)(\delta_\omega)
= v(\omega) = 1.$$
Similarly $V(\delta_\omega) = -1$ if $\omega \in Q$.
It follows that $f = 0$ on $H$, and so $f = 0$
on $S^2$.

To see that one may obtain many different cones on $Z''$ restricting to
the trivial cone on $Z$, simply choose a different partition $P', Q'$
of $H$.  It is easy to argue that the associated cone on $Z''$ is necessarily 
different to ${\mathfrak c}_u$ above.

The existence of multiple orderings on $E = Z''$ restricting to
the same maximal cone is very similar.  Suppose that 
${\mathfrak c} = {\mathfrak c}_U$ is a maximal cone on 
$Z$ corresponding
to a maximal open set $U$ as in Corollary
\ref{ntoc},
and suppose that $v$ is the associated maximal open tripotent in ${\mathcal Z}(E)$.
Thus $p = v^2$ is open; let $q = 1-p$ be the complementary 
projection.    It is easy to see that
$E q$ contains many linearly independent selfadjoint
tripotents $w$ with $w^2 = q$ (if $E^2 q$ was 
one dimensional it is easy to argue that $F_C^2$ has codimension
1, that is it is a maximal ideal in $Z^2$.  
Here $F_C$ is as in \ref{trsidlb2}.
This forces $S^2 \setminus (U \cup (-U))
= \{ \zeta, - \zeta \}$ for some $\zeta \in S^2$, which
is absurd).  Any such $w$ gives rise as before
to a selfadjoint unitary $u \geq v
\in E$, and hence to a maximal
natural cone ${\mathfrak d}_u$ on $Z$ (see Theorem 
\ref{order_char}).  We claim that 
${\mathfrak c}_u \cap Z = {\mathfrak c}_U$.  
Indeed if $x \in Z$ then $x \in {\mathfrak c}_u$ if and only if 
$u x \geq 0$.  Similarly, if $x \in Z$ then
$x \in {\mathfrak c}_U$ 
if and only if $v x \geq 0$ and $p x = x$.
The latter condition implies that $u x = u v v x = v x \geq 0$.
Thus $ {\mathfrak c}_U \subset {\mathfrak c}_u \cap Z$.
 For the other direction, note that if 
$x \in Z$ and if $u x \geq 0$,
then as before  
$(vx)(\varphi) = x(\varphi) \geq 0$ for $\varphi \in U$.  Thus
$x \leq 0$ on $-U$.  It follows from 
Corollary \ref{ntoc} that $x \in {\mathfrak c}_U$. 
 
\vspace{3 mm}

\section{Causal structures}

\newcommand{\Aut}{\operatorname{Aut}}
\newcommand{\aut}{\operatorname{\mathfrak{aut}}}
\newcommand{\fk}{\mathfrak k}
\newcommand{\fp}{\mathfrak p}

The above results, especially those of section 4, apply to a
classification
of those causal structures on certain infinite dimensional manifolds
$M$ that
come from a quantization of the points of $M$, i.e.\ from an embedding of $M$
into a space of bounded operators on some Hilbert space $H$. The finite
dimensional situation is apparently well understood (see e.g.\
\cite{CS}). We will go into
details elsewhere, and only sketch some of the main points here.

A (local) causal structure on a manifold $M$ consists of a field of
(regular)
cones $C_m$ contained in the tangent space $T_m(M)$. This definition is
motivated by general relativity: A path $\gamma:[a,b]\to M$ is causal iff
$\gamma'(t)\in C_{\gamma(t)}$. All curves in space-time that come from
existing particles do have this property, where the cones $C_m$ are the
forward light cones provided by the Lorentzian structure on $M$.

Let $E$ be a complex vector space. A bounded symmetric domain $D\sub E$ is a
connected open subset with the property that for each point $x\in D$
there is a
biholomorphic involution $s_x:D\to D$ with fixed point $x$. Denote by
$\Aut D$
the group of all biholomorphic automorphisms of $D$. Then $\Aut D$ is a
(real) Lie
group which acts homogeneously on $D$. The tangent space $T_x(D)$ of $D$
at a
point $x\in D$ can be identified with $E$. Fix $x\in D$, as well as a
regular cone
$C_x$ in $T_x(D)$, invariant under the action of the isotropy subgroup of
$\Aut D$ at $x$.  Let $C_{\overline{x}}=d\Phi(x)(C_x)$ for any
$\overline{x}\in D$ and $\Phi\in\Aut D$ with $\Phi(x)=x'$. Then the emerging
family of cones provides $D$ with a causal structure which actually is
invariant under the action of $\Aut D$.  Furthermore, any involution (`real
form') on $D$ gives rise to further symmetric spaces, and suitably
chosen cones
give rise to causal structures on the latter.

The connection to the results presented in the previous sections consists in
the following: The Lie algebra $\aut D$ of $\Aut D$ can be identified with
certain vector fields $X:D\to E$. We furthermore have a splitting
$$
\aut D=\fk\oplus\fp,
$$
with $\fk= \{ X\in\aut D : X\theta=-X \}$, $\fp= \{ X\in\aut
D : X\theta=X \}$ and
$\theta: D\to D$ defined by $\theta(x)=-x$. In a large number of cases, the
vector space $\fp$ carries the structure of a TRO in a natural way, and
this is
where the results of the present paper may be applied.

  \end{document}